\documentclass[12pt]{article}
\usepackage[margin=2cm]{geometry}
\usepackage{VCMFREsty,subfiles,titling}
\usepackage{apacite}
\usepackage{etoolbox}
\makeatletter
\patchcmd{\nocite}{\ifx\@onlypreamble\document}{\iftrue}{}{}
\makeatother

\pdfoutput=1
\makeatletter
\def\writefile[2]{}

\begin{document}
	
	\title{High-Dimensional Varying Coefficient Models with Functional Random Effects}
	
	\author{Michael Law \and Ya\hspace{-.1em}'\hspace{-.1em}acov Ritov \thanks{Supported in part by NSF Grants DMS-1646108, DMS-1712962, and DMS-2113364.}}
	\date{
		University of Michigan\\
		\today
	}
	
	\renewcommand\footnotemark{}
	
	\maketitle
	
	\begin{abstract}
		We consider a sparse high-dimensional varying coefficients model with random effects, a flexible linear model allowing covariates and coefficients to have a functional dependence with time. For each individual, we observe discretely sampled responses and covariates as a function of time as well as time invariant covariates. Under sampling times that are either fixed and common or random and independent amongst individuals, we propose a projection procedure for the empirical estimation of all varying coefficients. We extend this estimator to construct confidence bands for a fixed number of varying coefficients.
	\end{abstract}
	
	\section{Introduction}\label{sectionintroduction}
	
	Consider the following varying coefficients model with functional random effects given by
	\begin{align}\label{modely}
	y_{i}(t_{i,j}) = x_{i}^\T \beta(t_{i,j}) + [z_{i}(t_{i,j})]^\T \gamma(t_{i,j}) + \xi_{i}(t_{i,j}) + \epsilon_{i}(t_{i,j})
	\end{align}
	for $i = 1, \dots, n$ and $j = 1, \dots, m_{i}$.  Here, $x_{i} \in \R^{p}$ is a vector of time invariant covariates, $z_{i}(\cdot) : (0, 1) \to \R^{q}$ is a vector function representing the time varying covariates, and $\beta(\cdot) : (0, 1) \to \R^{p}$ and $\gamma(\cdot) : (0,1) \to \R^{q}$ are time varying coefficients.  Moreover, $\xi_{i}(\cdot) : (0,1) \to \R$ is a continuous time mean zero stochastic process representing the individual random effect and $\epsilon_{i}(\cdot) : (0,1) \to \R$ is an independent error.  Finally, the values $t_{i,j} \in (0, 1)$ are the sampling times.
	
	The model in \eqref{modely} is useful for longitudinal data, where for the $i$th individual, we record $m_{i}$ observations over time.  Traditionally, varying coefficients models for longitudinal data consider errors that are mean zero stochastic process with unknown covariance structure, such as \cite{hoover1998}, whereas we partition the error into continuous time individual random effects, $\xi_{i}(t_{i,j})$, and independent errors, $\epsilon_{i}(t_{i,j})$.  Such a partitioning is reasonable whenever the mean function, $\xi_{i}(\cdot)$, for each individual is smooth.  The model in equation \eqref{modely} admits many special cases both in the low-dimensional setting and the ultra high-dimensional setting; we list two such examples below:
	
	\begin{enumerate}
		\item Let $p = 1$, $q = 0$, and $x_{i} = 1$ for all $i = 1, \dots, n$, yielding the model
		\begin{align}\label{modelcy}
		y_{i}(t_{i,j}) = \beta(t_{i,j}) + \xi_{i}(t_{i,j}) + \epsilon_{i}(t_{i,j}).
		\end{align}
		Then, the problem of estimating $\beta(\cdot)$ is equivalent to the mean function estimation problem from \cite{cai2011}.
		
		\item Let $p > n$, $q = 0$, and $m_{i} = 1$ for all $i = 1, \dots, n$, yielding the model
		\begin{align}\label{modelkp}
		y_{i}(t_{i}) = x_{i}^\T \beta(t_{i}) + \epsilon_{i}(t_{i}).
		\end{align}
		Here, in the model, since we only have a single observation per individual, we slightly abuse notation and write $\epsilon_{i}(t_{i})$ to denote $\xi_{i}(t_{i}) + \epsilon_{i}(t_{i})$.  If the function $\beta(\cdot)$ is sparse, then this corresponds to the sparse varying coefficients model considered by \cite{klopp2015}.
	\end{enumerate}
	
	The goal of the present paper is estimation and inference for the varying coefficients $\beta(\cdot)$ and $\gamma(\cdot)$.  Much of the extant literature on high-dimensional varying coefficient models focus on estimation and variable selection, such as \cite{wei2011}, \cite{lian2012}, \cite{xue2012}, \cite{klopp2015}, and \cite{lee2016}.  The problem of inference is less well understood in the high-dimensional setting.  To the best of our knowledge, the only paper exploring this problem is \cite{chen2018}, who only consider local hypothesis testing.  However, all of the current works assume that the time varying covariates, $z_{i}(\cdot)$, are independent of the functional random effects, $\xi_{i}(\cdot)$.  In practice, this assumption is not necessarily satisfied.  To motivate this, consider modeling the average height of five year old boys in different countries over time.  For each country, our time varying covariates consist of variables, such as Human Development Index (HDI), health expenditure, and urbanization rate, that are likely to be correlated with the functional random effect of country, which encapsulates, among other things, environmental factors.  This example is revisited in Section \ref{sectionnumerics}.
	Moreover, most of the current works in high-dimensions assume either independent errors (cf. \cite{wei2011}, \cite{xue2012}, and \cite{klopp2015}) or a single individual, ie. $n = 1$ but $m \to \infty$ (cf. \cite{lee2016}).  The work most similar to ours is \cite{bai2019}, who assume $m_{i}$ observations for individual, $i = 1, \dots, n$.  However, they only consider random and independent sampling times with correlated Gaussian errors in a Bayesian paradigm.	
	
	Thus, our contribution to high-dimensional varying cofficient models is fourfold:  (i) estimation of $\beta(\cdot)$ in the presence of $\xi_{i}(\cdot)$, (ii) estimation of $\gamma(\cdot)$ under dependence of $z_{i}(\cdot)$ and $\xi_{i}(\cdot)$, (iii) estimation of both $\beta(\cdot)$ and $\gamma(\cdot)$ under both random and independent or fixed and common sampling times, and (iv) construction of confidence bands for single components of $\beta(\cdot)$ and $\gamma(\cdot)$ with and without dependence.
	
	To these ends, we propose a framework for estimation and inference in both the low-dimensional and the high-dimensional setting.  Our estimators utilize orthogonal series to leverage recent developments in the high-dimensional linear models literature.  We revisit the two examples in equations \eqref{modelcy} and \eqref{modelkp} later when analyzing the convergence rate of our estimators.
	
	\subsection{Organization of the Paper}
	We end this section with a description of the notation that is used in the remainder of the paper.  Since our estimator has two-stages, we consider each stage separately.  In Section \ref{sectionp=0}, we consider the special case when there are no time invariant covariates ($p = 0$) while in Section \ref{sectionq=0}, we consider the setting where there are no time varying covariates ($q=0$).  Then, we combine these results together into the two-stage estimator in Section \ref{sectioncomb}.  In Section \ref{sectioninference}, we consider the problem of constructing confidence bands for a fixed varying coefficient.  Finally, Section \ref{sectionnumerics} provides the simulations and Section \ref{sectionrda} presents an analysis of the height data mentioned above.  For ease of presentation, we defer all of the proofs to the Supplement.  In addition, the Supplement contains the assumptions for Section \ref{sectionq=0}, the additional results of the simulations, and an analysis of yeast cell cycle data.
	
	\subsection{Notation and Definitions}\label{sectionnotation}
	
	Throughout, all of our variables have a dependence on $n$, but, when it does not cause confusion, we suppress this dependence.  Since our interest is mainly asymptotic, we adopt a random design regression framework embedded in a triangular array.  The vector of covariates $x_{i}$ is drawn from a distribution that does not depend on time whilst the vector of covariates $z_{i}(t_{i,j})$ is drawn from a distribution conditioned on the sampling times.  Furthermore, the $\xi_{i}(\cdot)$s are independent realizations of smooth mean zero stochastic processes and $\epsilon_{i}(\cdot)$s are random errors.  Then, we write $\e$ to denote the expectation with respect to the joint probability measure of $(\xi_{i}(\cdot))_{i=1}^{n}$ and $(\epsilon_{i}(\cdot))_{i=1}^{n}$.
	
	Regarding the sampling times, there are two commonly used paradigms: (i) independent random sampling times and (ii) common fixed sampling times.  In the first setting, we assume that the time points are all independently sampled from a distribution $f$ on $(0, 1)$, which is bounded away from zero and infinity.  In this setting, $\t$ denotes the expectation with respect to the sampling times $t_{i,j}$.  On the other hand, for the common sampling times, we assume that $m_{i} = m$ and $t_{i,j} = j / m$ for all $i = 1, \dots, n$ and $j = 1, \dots, m$, viewing the sampling times as deterministic.  In either case, for a fixed value of $i = 1, \dots, n$, we write $(t_{i,(j)})_{j=1}^{m_{i}}$ to denote the order statistics for $(t_{i,j})_{j=1}^{m_{i}}$.
	
	As mentioned in the Introduction, we consider an orthogonal series estimator.  For technical convenience, we use the trigonometric basis since the functions are uniformly bounded by $\sqrt{2}$.  There are many definitions of the trigonometric basis, but we use the following definition as in \cite{tsybakov2008}.
	
	\begin{definition}\label{definitiontrigbasis}
		For $t \in (0,1)$, the trigonometric basis functions, denoted by $(\phi_{k}(\cdot))_{k=1}^{\infty}$, are given by
		\begin{align*}
		\phi_{k}(t) \defined
		\begin{cases}
		1, & k = 1\\
		\sqrt{2} \cos(\pi k t), & k = 2, 4, \dots \\
		\sqrt{2} \sin(\pi (k - 1) t), & k = 3, 5, \dots.
		\end{cases}
		\end{align*}
	\end{definition}
	Occasionally, it is useful to view these as functions on the complex plane; we write $\Im$ to denote the imaginary unit.  Later, to simplify notation, we assume that both the varying coefficients and random effects are in the same periodic Sobolev class with smoothness $\alpha$, denoted by $\Wper(\alpha, R)$ (for example, see Definition 1.11 of \cite{tsybakov2008}).  Then, the functions $\beta(\cdot)$, $\gamma(\cdot)$, and $\xi_{i}(\cdot)$ admit an expansion over the trigonometric basis.  Let $\beth_{k} \in \R^{p}$ (Hebrew letter Bet), $\gimel_{k} \in \R^{q}$ (Hebrew letter Gimel), and $\daleth_{i,k} \in \R$ (Hebrew letter Samek) for $k = 1, 2, \dots$ denote the Fourier coefficients of $\beta(\cdot)$, $\gamma(\cdot)$, and $\xi_{i}(\cdot)$ for  $i = 1, \dots, n$.  Then, we may write
	\begin{align}\label{equationbasisexpansion}
	\beta(\cdot) = \sum_{k=1}^{\infty} \beth_{k} \phi_{k}(\cdot),
	&& \gamma(\cdot) = \sum_{k=1}^{\infty} \gimel_{k} \phi_{k}(\cdot),
	&& \xi_{i}(\cdot) = \sum_{k=1}^{\infty} \daleth_{i,k} \phi_{k}(\cdot).
	\end{align}
	Let $\betahat(\cdot)$, with expansion
	\begin{align*}
	\betahat(t) = \sum_{k=1}^{\infty} \bethhat_{k} \phi_{k}(t),
	\end{align*}
	denote an arbitrary estimator for $\beta(\cdot)$.  To evaluate $\betahat$, we consider either integrated squared error (ISE) defined by
	\begin{align*}
	\ise(\betahat) \defined \int_{0}^{1} \left( \betahat(t) - \beta(t) \right)^{2} dt
	\end{align*}
	or mean integrated squared error (MISE), where $\mise (\betahat) \defined \t \e ( \ise (\betahat) )$.  We use MISE for the low-dimensional estimators and bound ISE for high-dimensional estimators with high probability.  By Parseval's Theorem, it follows that integrated squared error is equivalent to
	\begin{align}\label{equationmiseparseval}
	\ise(\betahat) = \sum_{k=1}^{\infty} \left\Vert \bethhat_{k} - \beth_{k} \right\Vert_{2}^{2}.
	\end{align}
	It follows from Proposition 1.14 of \cite{tsybakov2008} that
	\begin{align}\label{equationtsybakov}
	\sum_{k=\kbeta + 1}^{\infty} \left\Vert \beth_{k} \right\Vert_{2}^{2} = \O(\sbeta \kbeta^{-2\alpha}).
	\end{align}
	Therefore, it suffices to estimate the Fourier coefficients up to a truncation level $\kbeta$ to balance the bias-variance tradeoff.  Similarly, for $\gamma(\cdot)$, we find a truncation level $\kgamma$.  Thus, we may then define the low and high frequency components of the varying coefficient functions $\beta(\cdot)$ and $\gamma(\cdot)$ as
	\begin{align*}
	\betalower(\cdot) \defined \sum_{k=1}^{\kbeta} \beth_{k} \phi_{k}(\cdot), &&& \betaupper(\cdot) \defined \sum_{k=\kbeta + 1}^{\infty} \beth_{k} \phi_{k}(\cdot), \\
	\gammalower(\cdot) \defined \sum_{k=1}^{\kgamma} \gimel_{k} \phi_{k}(\cdot), &&& \gammaupper(\cdot) \defined \sum_{k=\kgamma + 1}^{\infty} \gimel_{k} \phi_{k}(\cdot).
	\end{align*}
	
	Similar to other works in high-dimensional statistics, we assume that our errors are sub-Gaussian, for which we use the following definition.
	\begin{definition}\label{definitionsubgaussian}
		A random vector $\eta \in \R^n$ is said to be \emph{sub-Gaussian} with parameter $\sgparam^{2}$, denoted $\eta \sim \sg_n \left( \sgparam^2 \right)$, if
		\begin{align*}
		\e \exp\left( \lambda^\T \eta \right) \leq \exp\left( \frac{\sgparam^2\left\Vert \lambda \right\Vert^2}{2} \right)
		\end{align*}
		for all $\lambda \in \R^n$.		
	\end{definition}
	
	Like other works in high-dimensional statistics, sparsity plays a crucial role.  For simplicity, we assume the setting of strong sparsity, whereby both $\beta(\cdot)$ and $\gamma(\cdot)$ have $\sbeta$ and $\sgamma$ components that are nonzero.  When considering the inferential problem, we need another notion of sparsity from \cite{vandegeer2014}.  Let $\Sigma$ to denote the population covariance matrix of $x_{i}$, $\Theta$ the inverse of $\Sigma$, and $\stheta = \max_{j = 1, \dots, p}|\{k \neq j : \Theta_{j,k} \neq 0 \}|$.  Thus, $\stheta$ is the maximal sparsity when regressing a component of $x_{i}$ against the remaining $x_{i}$'s.
	
	\section{Estimation with No Time Invariant Covariates}\label{sectionp=0}
	
	In this section, we assume that $p = 0$.  That is, the model we consider is
	\begin{align}\label{modelp=0}
	y_{i}(t_{i,j}) = [z_{i}(t_{i,j})]^\T \gamma(t_{i,j}) + \xi_{i}(t_{i,j}) + \epsilon_{i}(t_{i,j}).
	\end{align}
	Since the processes $z_{i}(\cdot)$ and $\xi_{i}(\cdot)$ may have arbitrary dependence, to remove the effect of $\xi_{i}(\cdot)$ from the model, we difference the observations that are sufficiently close.  As the function $\xi_{i}(\cdot)$ is assumed to be smooth, the value of $\xi_{i}(t)$ is approximately constant in a small neighborhood of $t$.  Let $h > 0$ be a bandwidth tuning parameter.  For simplicity, we temporarily assume that $m_{i}$ is even for each $i = 1, \dots, n$.   Then, we may define the set $\Aset_{h}$ as
	\begin{align*}
	\Aset_{h} \defined \left\{ (i,j) : 1 \leq i \leq n, j \in \left\{ 1, 3, \dots, m_{i} - 1 \right\}, t_{i, (j+1)} - t_{i,(j)} < h \right\}.
	\end{align*}
	Let $N \defined N_{h} = | \Aset_{h} |$.  For $(i,j) \in \Aset_{h}$, define the differenced observations $\ydiff_{i,j}$ as
	\begin{align*}
	\ydiff_{i,j} \defined \phant & y_{i}(t_{i,(j+1)}) - y_{i}(t_{i,(j)}) \\
	=\phant & [z_{i}(t_{i,(j+1)})]^\T \gamma(t_{i,(j+1)}) - [z_{i}(t_{i,(j)})]^\T \gamma(t_{i,(j)}) + \xi_{i}(t_{i,(j+1)}) - \xi_{i}(t_{i,(j)}) \\
	&+ \epsilon_{i}(t_{i,(j+1)}) - \epsilon_{i}(t_{i,(j)}) \\
	=\phant & \sum_{k=1}^{\kgamma} \underbrace{[\phi_{k}(t_{i,(j+1)}) z_{i}(t_{i,(j + 1)}) - \phi_{k}(t_{i,(j)}) z_{i}(t_{i,(j)})]^\T}_{\psi_{i,j,k}^\T} \gimel_{k} + \underbrace{\epsilon_{i}(t_{i,(j+1)}) - \epsilon_{i}(t_{i,(j)})}_{\eta_{i,j}} \\
	&+ \underbrace{\sum_{k=\kgamma +1}^{\infty} [\phi_{k}(t_{i,(j+1)}) z_{i}(t_{i,(j + 1)}) - \phi_{k}(t_{i,(j)}) z_{i}(t_{i,(j)})]^\T \gimel_{k} + \xi_{i}(t_{i,(j+1)}) - \xi_{i}(t_{i,(j)})}_{\remgamma_{i,j}} \\
	= \phant & \psi_{i,j}^\T \gimelvec + \eta_{i,j} + \remgamma_{i,j},
	\end{align*}
	where $\psi_{i,j} = (\psi_{i,j,1}^\T , \dots, \psi_{i,j,\kgamma}^\T)^\T$ and $\gimelvec = (\gimelvec_{1}^\T, \dots, \gimelvec_{\kgamma}^\T)^\T$.  In matrix notation, we write this model as
	\begin{align}\label{modelplmgimelvec}
	\Ydiff = \Psi \gimelvec + \etavec + \remgammavec.
	\end{align}
	This is a sparse high-dimensional partially linear model with uncorrelated errors, for which there are many proposals for estimating $\gimelvec$.  Commonly, in high-dimensional nonparametric models, a version of group lasso (cf. \cite{yuan2006}) is used to select relevant functions after a basis expansion, such as the SpAM estimator of \cite{ravikumar2009} or the block lasso estimator of \cite{klopp2015}.  While such approaches lead to more interpretable estimators, we estimate $\gimelvec$ by the classical lasso of \cite{tibshirani1996} but note that our approach generalizes to using the group lasso.  We use the classical lasso to motivate the inferential procedure in Section \ref{sectioninference}, which is based on a version of the de-biased lasso estimator.  Therefore, we estimate $\gimelvec$ in the low-dimensional case by
	\begin{align*}
	\gimelvechatld \defined \left( \Psi^\T \Psi \right)^{-1} \Psi^\T \Ydiff
	\end{align*}
	and
	\begin{align*}
	\gimelvechathd \defined \argmin_{\gimelvec \in \R^{q\kgamma}} N^{-1} \left\Vert \Ydiff - \Psi \gimelvec \right\Vert_{2}^{2} + \lambda \left\Vert \gimelvec \right\Vert_{1}
	\end{align*}
	in the high-dimensional setting for a suitable tuning parameter $\lambda > 0$.  By identifying the vectors $\gimelhatld$ and $\gimelhathd$ as $\kgamma$ vectors in $\R^{q}$, the estimators for $\gamma(\cdot)$ are given by
	\begin{align*}
	\gammahatld(\cdot) \defined \sum_{k=1}^{\infty} \kgamma \gimelhatld_{k} \phi_{k}(\cdot), &&
	\gammahathd(\cdot) \defined \sum_{k=1}^{\infty} \kgamma \gimelhathd_{k} \phi_{k}(\cdot).
	\end{align*}

	\begin{remark}
		In the above formulation, we pair the observations in $\Aset_{h}$ to ensure that the resultant errors $\etavec$ are independent.  However, this reduces the number of observations that we have to estimate $\gimelvec$.  To circumvent this problem, we may alternatively consider the set
		\begin{align*}
		\Bset_{h} \defined \left\{ (i,j) : 1\leq i \leq n, 1 \leq j \leq m_{i} - 1, t_{i,(j + 1)} - t_{i,(j)} < h \right\}.
		\end{align*}
		Using the set $\Bset_{h}$, we may likewise form the model given in equation \eqref{modelplmgimelvec}, where the resultant partially linear model has, by construction, correlated errors with known correlation structure.  Hence we may find a matrix $B$ such that
		\begin{align*}
		B \Ydiff = B\Psi \gimelvec + B\etavec + B\remgammavec,
		\end{align*}
		where $B\etavec$ is uncorrelated.  For simplicity, we consider only the set $\Aset_{h}$, but in practice, we recommend adjusting using $\Bset_{h}$ when the sampling times are random and independent but $\Aset_{h}$ when the sampling times are fixed and common (see Table \ref{tablegammaq500} in Section \ref{sectionsupplementsimulations} of the Supplement).
	\end{remark}
	
	\subsection{Sample Size}\label{sectionp=0samplesize}
	
	In this subsection, we consider the expected number of observations after differencing under a few asymptotic regimes for $n$, $m$, and $h$.  This leads us to the following proposition.
	
	\begin{proposition}\label{propositionsamplesize}
		Suppose the sampling times $t_{i,j} \iid f$ for a density $f$ on $(0, 1)$ bounded away from zero and infinity and $m_{i} = m > 0$ for all $i = 1, \dots, n$.  Let $\Ntilde_{h} = |\Bset_{h}|$.
		
		\begin{enumerate}
			\item If $m = \O(1)$ and $n \to \infty$, then $\t \Ntilde_{h} \asymp n h$.
			
			\item If $m \to \infty$ and $mh \ll 1$, then $\t \Ntilde_{h} \asymp n m^{2} h$ and $\t (N_{h} + 1)^{-1} \asymp (nm^{2}h)^{-1}$.
			
			\item If $m \to \infty$ and $mh \gg 1$, then $\t \Ntilde_{h} \asymp nm$.  If, in addition, $mh - \log(mn) \to \infty$, then $\p (\Ntilde_{h} = mn) \to 1$.
		\end{enumerate}
	\end{proposition}
	
	\begin{remark}
		It is easy to see that $| \Aset_{h} | \asymp | \Bset_{h} |$.
	\end{remark}
	
	\subsection{Assumptions}
	
	The following assumptions are used when $p = 0$.
	
	\begin{enumerate}[label=(A\arabic*)]
		\item \label{assumptionp=0designld} If $\sgamma = q < N$, then the matrix $\Psi$ satisfies $\tr [( \Psi^\T \Psi )^{-1} ] = \O( \sgamma \kgamma / N)$ and $\Vert (\Psi^\T \Psi)^{-1} \Vert_{2} = \Op(N^{-1})$.
		
		\item \label{assumptionp=0designscaling} The columns of the matrix $\Psi$ have squared norms that are uniformly $\Op(N)$.
		
		\item \label{assumptionp=0psicc} The design matrix $\Psi$ satisfies the compatibility condition with compatibility constant $\ccpsi > 0$.
		
		\item \label{assumptionp=0psiarev} The design matrix $\Psi$ satisfies the adaptive restricted eigenvalue condition with constant $\arevpsi > 0$.
		
		\item \label{assumptionp=0errors} The errors $\epsilon_{i}(t_{i,j})$ are independent and identically distributed with mean zero and variance $\sigmaepsilon$.  Moreover, the errors are independent of the sampling times.
		
		\item \label{assumptionp=0sgerrors} The errors $\epsilon_{i}(t_{i,j}) \iid \sg(\sgparam_{\epsilon}^{2})$.  Moreover, the errors are independent of the sampling times.
		
		\item \label{assumptionp=0lipschitz} The functional random effects $(\xi_{i}(\cdot))_{i=1}^{n}$ are uniformly Lipschitz with constant $\lip$.
		
		\item \label{assumptionp=0gamma} Each coordinate of the coefficient $\gamma(\cdot)$ satisfies $\gamma_{k}(t) \in \Wper(\alpha, R)$ for some constant $\alpha \geq 2$ and $R > 0$ with $\sgamma$ coordinates nonzero.  Moreover,
		\begin{align*}
		\t \left( \sum_{k= \kgamma + 1}^{\infty}[z_{i} (t)]^\T \gimel_{k} \phi_{k}(t) \right)^2 = \O( \sgamma \kgamma^{-2\alpha}).
		\end{align*}
	\end{enumerate}
	
	\begin{remark}
		Assumptions \ref{assumptionp=0designld}, \ref{assumptionp=0designscaling}, and \ref{assumptionp=0errors} are standard scaling assumptions for the design in the low-dimensional setting.
		
		Assumptions \ref{assumptionp=0psicc} and \ref{assumptionp=0psiarev} are both compatibility conditions on the design matrix, with \ref{assumptionp=0psiarev} implying \ref{assumptionp=0psicc}.  In Theorem \ref{theoremp=0hd} below, we use \ref{assumptionp=0psicc} to obtain slow rates on ISE while \ref{assumptionp=0psiarev} yields a fast rate on ISE.  For a more detailed discussion on assumptions \ref{assumptionp=0psicc} and \ref{assumptionp=0psiarev}, including definitions, we refer the reader to Section 6.2 of \cite{buhlmann2011} and Section 4 of \cite{bickel2009} respectively.
		
		Next, assumption \ref{assumptionp=0sgerrors} is standard in the high-dimensional linear models literature and \ref{assumptionp=0lipschitz} is reasonable whenever the underlying mean function for each individual is smooth.  Further, \ref{assumptionp=0lipschitz} is implied whenever $(\xi_{i}(\cdot))_{i=1}^{n} \subseteq \Wper(\alpha, R)$.
		
		The first half of assumption \ref{assumptionp=0gamma} is standard in the literature on nonparametric regression while the second half ensures that the varying coefficients can be well approximated by a few basis functions.  In Example \ref{exampletimerandomcovariates} below, we consider an instance where the second half is satisfied.
	\end{remark}
	
	\begin{example}[Example for Assumption \ref{assumptionp=0gamma}:  Time random covariates]\label{exampletimerandomcovariates}
		Suppose that the distribution of $z_{i}(t_{i,j})$ does not depend on $t_{i,j}$.  That is, assume that $\dist{z_{i}(t_{i,j}) | t_{i,j}} = \dist{z_{i}(t_{i,k}) | t_{i,k}}$ for every $j, k = 1, \dots, m_{i}$ with $j \neq k$.  Then,
		\begin{align*}
		\t \left( \sum_{k= \kgamma + 1}^{\infty}[z_{i} (t)]^\T \gimel_{k} \phi_{k}(t) \right)^2
		&= \sum_{k=\kgamma+1}^{\infty} \t \left( [z_{i}(t)]^\T \gimel_{k} \right)^2 \\
		&= \sum_{k=\kgamma+1}^{\infty} \gimel_{k}^\T \t \left( [z_{i}(t)] [z_{i}(t)]^\T \right) \gimel_{k} \\
		&\leq \left\Vert \t \left( [z_{i}(t)] [z_{i}(t)]^\T \right) \right\Vert_{2}^{2} \sum_{k=\kgamma+1}^{\infty} \left\Vert \gimel_{k} \right\Vert_{2}^{2} \\
		&= \O\left( \sgamma \kgamma^{-2\alpha} \right).
		\end{align*}
	\end{example}
	
	\subsection{Main Results}\label{sectionp=0estimation}
	
	We start by stating a result for the low-dimensional setting.
	
	\begin{proposition}\label{theoremp=0ld}
		Consider the model given in equation \eqref{modelplmgimelvec}.  Assume \ref{assumptionp=0designld}, \ref{assumptionp=0errors}, \ref{assumptionp=0lipschitz}, and \ref{assumptionp=0gamma}.  Then
		\begin{align*}
		\e \t \left\Vert \gimelvechatld - \gimelvec \right\Vert_{2}^{2}
		= \O \left( \sgamma \kgamma \t N^{-1} + \sgamma \kgamma^{-2\alpha} + \lip^{2} h^{2} \right).
		\end{align*}
	\end{proposition}
	
	\begin{remark}
		As noted in Section \ref{sectionnotation}, the MISE of $\gammahatld(\cdot)$ can be bounded by
		\begin{align*}
		\mise(\gammahatld) = \O \left( \sgamma \kgamma \t N^{-1} + \sgamma \kgamma^{-2\alpha} + \lip^{2} h^{2} \right).
		\end{align*}
		Choosing $\kgamma \asymp (\t N^{-1})^{-1 / (2\alpha + 1)}$ yields
		\begin{align*}
		\mise(\gammahatld) = \O \left( \sgamma (\t N^{-1})^{2\alpha / (2\alpha + 1)} + \lip^{2} h^{2} \right).
		\end{align*}
		The choice of $h$ is less straightforward as it depends on the asymptotic growth of $m$ relative to $n$, which can be seen from Proposition \ref{propositionsamplesize}.
	\end{remark}
	
	Now, turning our attention to the high-dimensional setting, we have the following result.
	
	\begin{theorem}\label{theoremp=0hd}
		Consider the model given in equation \eqref{modelplmgimelvec}.  Assume \ref{assumptionp=0designscaling}, \ref{assumptionp=0sgerrors}, \ref{assumptionp=0lipschitz}, and \ref{assumptionp=0gamma}.  For $t > 0$, let
		\begin{align}
		\lambda_{0} \defined 2 \sgparam_{\epsilon} \sqrt{N^{-1} \max_{j=1, \dots, p} \left\Vert \Psi_{j} \right\Vert_{2}^{2}} \sqrt{ \frac{t^{2} + 2\log(q\kgamma)}{N}}.
		\end{align}
		Suppose $\lambda \geq 2\lambda_{0}$.
		
		\begin{enumerate}
			\item If, in addition, \ref{assumptionp=0psicc} holds, then with probability at least $1 - 2\exp(-t^{2} / 2)$,
			\begin{align*}
			2 N^{-1} &\left\Vert \Psi \gimelvechathd - \Psi \gimel - \remgammavec \right\Vert_{2}^{2}
			+ \lambda \left\Vert \gimelvechathd - \gimelvec \right\Vert_{1} \\
			&\leq 6 N^{-1} \left\Vert \remgammavec \right\Vert_{2}^{2}
			+ 24 \ccpsi^{-2} \lambda^{2} \sgamma \kgamma.
			\end{align*}
			
			\item If, in addition to the above, \ref{assumptionp=0psiarev} holds, then with probability at least $1 - 2\exp( -t^{2} / 2)$,
			\begin{align*}
			\left\Vert \gimelvechathd - \gimelvec \right\Vert_{2}^{2}
			= \O \left( \lambda^{2} \sgamma \kgamma \left( \frac{\sgamma \kgamma^{-2\alpha} + \lip^{2}h^{2}}{\lambda^{2} \sgamma \kgamma} + \arevpsi^{-2} \right)^{2} \right).
			\end{align*}
		\end{enumerate}
	\end{theorem}
	
	\setcounter{corollaryp=0}{\value{theorem}}
	
	\begin{example}
		As a special case, consider the setting where $n = 1$ and $m = m_{1} \to \infty$.  Note that this problem is a generalization of the model in equation \eqref{modelkp} from the Introduction.  Under the additional assumption that $\xi_{i}(\cdot) \equiv 0$, obtaining $m$ observations from a single individual with time varying covariates is equivalent to obtaining a single observation from $m$ individuals.  Set $\lambda = 2\lambda_{0}$, $\kgamma \asymp (m / \log(q))^{1 / (2\alpha + 1)}$ and $h \asymp \sgamma (\log(q) / m)^{2\alpha / (2\alpha + 1)}$.  Since $q > \kgamma$, it follows that $\log(q) \leq \log(q\kgamma) \leq 2 \log(q)$.  Moreover, the choice of $h$ implies that $mh \gg 1$; thus, $N = m$ with high probability for $m$ sufficiently large by Proposition \ref{propositionsamplesize}.  Then, with probability at least $1 - 2\exp(-t^{2} / 2)$, it follows from Theorem \ref{theoremp=0hd} that
		\begin{align*}
		\left\Vert \gimelhathd - \gimelvec \right\Vert_{2}^{2}
		&= \O \left(
		\frac{\sgamma m}{\kgamma^{4\alpha + 1}\log(q \kgamma)}
		+ \frac{h^{2} m}{\sgamma \kgamma \log(q \kgamma)}
		+ \frac{\sgamma \kgamma \log(q \kgamma)}{m}
		\right)\\
		&= \O \left( \sgamma \left(\frac{\log(q)}{m}\right)^{2\alpha / (2\alpha + 1)} \right).
		\end{align*}
		From equation \eqref{equationtsybakov}, it follows that
		\begin{align*}
		\sum_{k=\kgamma + 1}^{\infty} \left\Vert \gimel_{k} \right\Vert_{2}^{2} = \O\left( \sgamma \kgamma^{-2\alpha} \right).
		\end{align*}
		Combining these facts yields the following bound on ISE with probability at least $1 - 2\exp(-t^{2} / 2)$.
		\begin{align*}
		\ise(\gammahathd)
		= \left\Vert \gimelhathd - \gimelvec \right\Vert_{2}^{2}
		+ \sum_{k=\kgamma + 1}^{\infty} \left\Vert \gimel_{k} \right\Vert_{2}^{2}
		= \O \left( \sgamma \left(\frac{\log(q)}{m}\right)^{2\alpha / (2\alpha + 1)} \right).
		\end{align*}
		From Theorem 1 of \cite{klopp2015}, assuming that the functions have uniform smoothness $\alpha$, then, up to the logarithmic factor, $\gammahathd(\cdot)$ attains the minimax rate.		
	\end{example}

	\section{Estimation with No Time Varying Covariates}\label{sectionq=0}
	
	In this section, we assume that $q = 0$.  That is, the model we consider is
	\begin{align}\label{modelq=0}
	y_{i}(t_{i,j}) = x_{i}^\T \beta(t_{i,j}) + \xi_{i}(t_{i,j}) + \epsilon_{i}(t_{i,j}).
	\end{align}
	Now, by directly substituting the expansions from equation \eqref{equationbasisexpansion}, it follows that
	\begin{align*}
	y_{i}(t_{i,j}) = \sum_{k=1}^{\infty} \left( x_{i}^\T \beth_{k} + \daleth_{i,k} \right) \phi_{k}(t_{i,j}) + \epsilon_{i}(t_{i,j}).
	\end{align*}
	The above factoring suggests that, to estimate the $k$'th Fourier coefficient $\beth_{k}$, we should look at the observations in the frequency domain as opposed to the time domain.  That is, we should form new observations in the frequency domain as
	\begin{align*}
	\yproj_{i,k}
	&\defined m_{i}^{-1} \sum_{j=1}^{m_{i}} y_{i}(t_{i,j}) \phi_{k}(t_{i,j}).
	\end{align*}
	By projecting the observations onto a fixed frequency given by $\phi_{k}$, the above is an approximate linear model.  Depending on whether the sampling times are viewed as fixed and common or random and independent, we partition the above model differently.
	
	In the setting where the sampling times are random and independent for each individual, we define
	\begin{align*}
	\zeta_{i, k} \defined m_{i}^{-1} \sum_{j=1}^{m_{i}} \left( x_{i}^\T \beta(t_{i,j}) + \xi_{i}(t_{i,j}) + \epsilon_{i}(t_{i,j}) \right) \phi_{k}(t_{i,j})- x_{i}^\T \beth_{k},
	\end{align*}
	yielding a linear model
	\begin{align}\label{modelfdind}
	\yproj_{i,k} = x_{i}^\T \beth_{k} + \zeta_{i,k}.
	\end{align}
	In matrix notation, we write
	\begin{align*}
	\Yproj_{k} = X \beth_{k} + \zetavec_{k}.
	\end{align*}	
	
	When the sampling times are fixed and common, we write
	\begin{align}\label{modelfdcom}
	\yproj_{i, k} = x_{i}^\T (\beth_{k} + \samekh_{k}) + \zeta_{i, k},
	\end{align}
	where
	\begin{align*}
	&\samekh_{k} \defined m^{-1} \sum_{\substack{l = 1 \\ l \neq k}}^{\infty} \beth_{l} \sum_{j = 1}^{m} \phi_{k}(t_{i,j}) \phi_{l}(t_{i,j}), \hspace{1em}(\text{Hebrew letter Dalet}) \\
	&\zeta_{i,k} \defined \daleth_{i, k} + m^{-1} \sum_{j = 1}^{m} \phi_{k}(t_{i,j}) \left( \epsilon_{i}(t_{i,j}) + \sum_{\substack{l = 1 \\ l \neq k}}^{\infty} \daleth_{i, l} \phi_{l}(t_{i,j}) \right).
	\end{align*}
	In matrix notation, this is expressed as
	\begin{align*}
	\Yproj_{k} = X (\beth_{k} + \samekh_{k}) + \zetavec_{k}.
	\end{align*}
	
	The main difference between the perspectives in equations \eqref{modelfdind} and \eqref{modelfdcom} is how we consider the inexact orthogonalization.  When the sampling times are random, the projection of $x_{i}^\T \beta(\cdot)$ over the $k$'th frequency is unbiased for $x_{i}^\T \beth_{k}$ with respect to the probability measure on $t_{i,j}$.  Conversely, since the common sampling times are deterministic, this inexact orthogonalization is viewed as bias; the vector $\samekh_{k}$ is the sum of the Fourier coefficients of the aliased frequencies.
	
	Regardless of the sampling times, in the low-dimensional setting, least-squares provides a convenient estimator for $\beth_{k}$, while in the high-dimensional setting, one may directly apply the lasso.  That is,
	\begin{align*}
	\bethhatld_{k} &\defined
	\begin{cases}
	\left( X^\T X \right)^{-1} X^\T \Yproj_{k}, & k = 1, \dots, \kbeta, \\
	\zerovec_{p},
	\phantom{\argmin_{\beth \in \R^{p}} n^{-1} \left\Vert \Yproj_{k} - X \beth \right\Vert_{2}^{2} + \lambda_{k} \left\Vert \beth \right\Vert_{1}}
	&  k > \kbeta,
	\end{cases}
	\\
	\bethhathd_{k} &\defined
	\begin{cases}
	\argmin_{\beth \in \R^{p}} n^{-1} \left\Vert \Yproj_{k} - X \beth \right\Vert_{2}^{2} + \lambda_{k} \left\Vert \beth \right\Vert_{1}, & k = 1, \dots, \kbeta, \\
	\zerovec_{p},
	\phantom{\argmin_{\beth \in \R^{p}} n^{-1} \left\Vert \Yproj_{k} - X \beth \right\Vert_{2}^{2} + \lambda_{k} \left\Vert \beth \right\Vert_{1}}
	& k > \kbeta.
	\end{cases}
	\end{align*}
	Like in Section \ref{sectionp=0}, this provides the estimators for $\beta(\cdot)$ as
	\begin{align*}
	\betahatld(\cdot) \defined \sum_{k=1}^{\infty} \bethhatld_{k} \phi_{k}(\cdot), &&
	\betahathd(\cdot) \defined \sum_{k=1}^{\infty} \bethhathd_{k} \phi_{k}(\cdot).
	\end{align*}
	
	\subsection{Assumptions}
	
	The assumptions are provided in Section \ref{sectionsupplementassumption} of the Supplement.
	
	There are two assumptions regarding the sampling times, \ref{assumptionq=0samplingtimes} and \ref{assumptionq=0samplingtimesdiscrete}.  As first pointed out by \cite{cai2011} in the context of mean function estimation, the rate of convergence is different between random and independent or fixed and common sampling times.  The first assumption, \ref{assumptionq=0samplingtimes}, considers the independent sampling times, with an additional convenience that the sampling times are uniform since the trigonometric basis is orthogonal with respect to this measure on $(0,1)$.  This assumption may be relaxed to allow for $t_{i,j} \iid f$ from some density $f$ bounded from zero and infinity by taking an appropriate change of measure.  The other assumption, \ref{assumptionq=0samplingtimesdiscrete}, considers the common sampling time setting.  Again, we make a simplifying assumption that the sampling times are on a uniformly spaced grid, though this assumption may similarly be relaxed.  These two settings are analyzed separately below.
	
	Assumptions \ref{assumptionq=0designscaling} -- \ref{assumptionq=0xarev} are analogous \ref{assumptionp=0designscaling} -- \ref{assumptionp=0psiarev}.
	
	Next, assumption \ref{assumptionq=0betaint} assumes that the signal is uniformly bounded by some function $g(n)$.  This assumption is the finite sample analogue of maintaining a bounded signal to noise ratio in a varying coefficients model.  We conflate the function $g(n)$ in assumptions \ref{assumptionq=0designscaling} and \ref{assumptionq=0betaint} since $g(n)$ may normally be taken to be a slowly varying function of $n$.  For example, if the design is sub-Gaussian, then $g(n) = \log(n)$ is sufficient.  If the design is bounded, then $g(n)$ may be a constant.  Further, we note that this assumption implies that
	\begin{align*}
	\sup_{i=1, \dots, n} \sup_{t \in (0,1)} \left( x_{i}^\T \beta(t) \right)^{2} dt = \O(g(n)).
	\end{align*}
	
	Assumption \ref{assumptionq=0xiint}  is a slightly stronger version of \ref{assumptionp=0lipschitz}.  It automatically implies that $\e \daleth_{i,k} = 0$ for all $k = 1, 2, \dots$, and $\sum_{k=\kbeta}^{\infty} \sigma_{\daleth,k}^{2} = \O(\kbeta^{-2\alpha})$ where $\sigma_{\daleth,k} = \var(\daleth_{i,k})$.  Finally, \ref{assumptionq=0xisup} is a technical requirement for the high-dimensional regime to ensure concentration of the resultant high-dimensional linear models after projection.
	
	\subsection{Main Results:  Independent Sampling Times}\label{sectionq=0independent}
	
	We start by considering the low-dimensional regime.
	\begin{proposition}\label{theoremq=0ld}
		Consider the model given in equation \eqref{modelq=0} with $\sbeta = p < n$.  Let $M = \diag((m_{1}, \dots, m_{n}))$.  Under assumptions \ref{assumptionq=0samplingtimes}, \ref{assumptionq=0betaint}, and \ref{assumptionq=0xiint}, the MISE is given by
		\begin{align*}
		\mise (\betahatld)
		=& \phant \tr\left( \left( X^\T X \right)^{-1} X^\T
		\left( \O\left( 1 \right) I_{n}
		+ \O(g(n) \kbeta)M^{-1} \right)
		X \left( X^\T X \right)^{-1} \right) \\
		&+ \O(\sbeta \kbeta^{-2\alpha}).
		\end{align*}
	\end{proposition}
	
	Then, $\kbeta$ is chosen to minimize the right hand side.  In general, the optimal value of $\kbeta$ is dependent on the specific sequence of designs.  We consider a special case where we can characterize the exact tradeoff between the number of unique individuals, $n$, and the number of samples per observation $m_{i}$.
	\begin{example}[Minimax Estimation of the Mean Function]\label{examplecy}
		Consider the model given in equation \eqref{modelcy} with independent uniform sampling times, which is reproduced below for convenience.
		\begin{align*}
		y_{i}(t_{i,j}) = \beta(t_{i,j}) + \xi_{i}(t_{i,j}) + \epsilon_{i}(t_{i,j}).
		\end{align*}
		Here, $\sbeta = p = 1$ with $x_{i} = 1$ for $i = 1, \dots, n$.  In this case, we may set $g(n) = 1$.  For now, we write $m = \left( \sum_{i=1}^{n} m_{i}^{-1} \right)^{-1}$ to denote the harmonic mean of the $( m_{i} )_{i=1}^{n}$.  Then,
		\begin{align*}
		\tr \left( \left( X^\T X \right)^{-1} X^\T I_{n}X \left( X^\T X \right)^{-1} \right) = \O\left(n^{-1}\right).
		\end{align*}
		Next, by a direct calculation,
		\begin{align*}
		\sum_{k=1}^{\kbeta} \tr \left( \left( X^\T X \right)^{-1} X^\T  M^{-1} X \left( X^\T X \right)^{-1} \right) = \O\left( \kbeta (mn)^{-1} \right)
		\end{align*}
		Thus, the risk from Proposition \ref{theoremq=0ld} can be simplified to
		\begin{align*}
		\mise(\betahatld) = \e \int_{0}^{1} \left\Vert \betahatld(t) - \beta(t) \right\Vert^{2} dt
		= \O\left(n^{-1} + \kbeta(mn)^{-1} + \kbeta^{-2\alpha}\right).
		\end{align*}
		This yields the optimal choice of $\kbeta \asymp (mn)^{1/(2\alpha + 1)}$.  Then, the risk is
		\begin{align*}
		\e \int_{0}^{1} \left\Vert \betahatld(t) - \beta(t) \right\Vert^{2} dt
		= \O\left(n^{-1} + (mn)^{-2\alpha/(2\alpha + 1)}\right),
		\end{align*}
		which coincides with the minimax rate from Theorem 3.1 of \cite{cai2011}.
	\end{example}
	
	The next result is the analogue of Proposition \ref{theoremq=0ld} for the high-dimensional setting.
	
	\begin{theorem}\label{theoremq=0hd}
		Consider the model given in equation \eqref{modelq=0}.  Assume \ref{assumptionp=0sgerrors}, \ref{assumptionq=0samplingtimes}, \ref{assumptionq=0designscaling}, \ref{assumptionq=0betaint}, and \ref{assumptionq=0xisup}.  For $k \leq \kbeta$ and $t > 0$, let
		\begin{align*}
		\lambda_{0, k} \defined \sqrt{\max_{j=1,\dots, p} n^{-1} \sum_{i=1}^{n} \sgparam_{\zeta, i, k}^{2} x_{i,j}^{2}} \sqrt{\frac{t^{2} + 2\log(p)}{n}},
		\end{align*}
		where $\sgparam_{\zeta,i,k}^{2} = \O(\sgparam_{\daleth, k}^{2} + g(n)m_{i}^{-1})$. Suppose $\lambda_{k} \geq 2 \lambda_{0, k}$.
		\begin{enumerate}
			\item If in addition \ref{assumptionq=0xcc} holds, then with probability at least $1 - 2\exp(-t^{2} / 2)$,
			\begin{align*}
			n^{-1} \left\Vert X \left( \bethhat_{k} - \beth_{k} \right) \right\Vert_{2}^{2} + \lambda_{k} \left\Vert \bethhat_{k} - \beth_{k} \right\Vert_{1} \leq 4 \lambda_{k}^{2} \sbeta / \ccx^{2}.
			\end{align*}
			
			\item If in addition \ref{assumptionq=0xarev} holds, then with probability at least $1 - 2\exp(-t^{2} / 2)$,
			\begin{align*}
			\left\Vert \bethhat_{k} - \beth_{k} \right\Vert_{2}^{2}
			=
			\O \left( \lambda_{k}^{2} \sbeta \arevx^{-4} \right).
			\end{align*}
		\end{enumerate}
	\end{theorem}
	\setcounter{corollaryq=0}{\value{theorem}}
	
	\begin{example}\label{exampleq=0hd}
		As a special case, consider the setting where $m_{i} = m$ for all $i = 1, \dots, n$.  Then,
		\begin{align*}
		\lambda_{0,k}^{2}
		= \O\left(
		\left( \sgparam_{\daleth, k}^{2} + g(n) m^{-1} \right) \frac{\log(p)}{n} \right).
		\end{align*}
		Thus, under assumptions \ref{assumptionq=0samplingtimes}, \ref{assumptionq=0designscaling}, \ref{assumptionq=0xarev}, \ref{assumptionq=0betaint}, and \ref{assumptionq=0xisup}, we have for $k = 1, \dots, \kbeta$ that
		\begin{align*}
		\left\Vert \bethhat_{k} - \beth_{k} \right\Vert_{2}^{2} = \O\left( \left( \sgparam_{\daleth, k}^{2} + g(n) m^{-1} \right) \frac{\sbeta \log(p)}{n} \right).
		\end{align*}
		Then, with probability at least $1 - 2\exp(-t^{2} / 2 + \log(\kbeta))$,
		\begin{align*}
		\ise(\betahathd) = \O\left( \left( 1 + \kbeta m^{-1} g(n) \right) \frac{\sbeta \log(p)}{n} + \sbeta \kbeta^{-2\alpha} \right).
		\end{align*}
		This yields an optimal choice of $\kbeta \asymp (nm / (g(n) \log(p))^{1/(2\alpha + 1)}$ with risk
		\begin{align*}
		\ise(\betahathd) = \O\left( \frac{\sbeta \log(p)}{n} + \sbeta \left( \frac{g(n) \log(p)}{nm} \right)^{2\alpha / (2\alpha + 1)} \right) .
		\end{align*}
		If we assume that $m = 1$, then this corresponds to the model in equation \eqref{modelkp}, which is a simplification of the model of \cite{klopp2015} under uniform smoothness.  Then, with probability at least $1 - 2\exp(-t^{2} / 2 + \log(\kbeta))$,
		\begin{align*}
		\ise( \betahathd ) = \O \left( \frac{\sbeta \log(p)}{n} + \sbeta \left( \frac{g(n) \log(p)}{n} \right)^{2\alpha / (2\alpha + 1)} \right),
		\end{align*}
		which, up to the slowly varying functions $g(n)$ and $\log(p)$, achieves the minimax rate from Theorem 1 of \cite{klopp2015}.
	\end{example}
	
	\setcounter{corollaryq=0}{\value{theorem}}

	\subsection{Main Results:  Common Sampling Times}\label{sectionq=0common}
	
	Again, we start by considering the low-dimensional estimator.
	
	\begin{proposition}\label{theoremq=0ldcommon}
		Consider the model given in equation \eqref{modelq=0} with $\sbeta = p < n$.  Under assumptions \ref{assumptionp=0errors}, \ref{assumptionq=0samplingtimesdiscrete}, \ref{assumptionq=0betaint}, and \ref{assumptionq=0xiint}, the MISE for $\kbeta \leq m -1$ is given by
		\begin{align*}
		\mise (\betahatld)
		= \O \left( \sbeta m^{-2\alpha} + (1 + \kbeta m^{-1}) \tr[(X^\T X)^{-1}] + \sbeta \kbeta^{-2\alpha} \right).
		\end{align*}
	\end{proposition}
	
	\begin{remark}
		In Proposition \ref{theoremq=0ldcommon}, if we make a scaling assumption analogous to \ref{assumptionp=0designscaling} with $\tr[(X^\T X)^{-1}] = \O(\sbeta / n)$, then by choosing $\kbeta \leq m -  1$ and $\kbeta \asymp m$, we may further obtain the bound
		\begin{align*}
		\mise(\betahatld) = \O \left( \sbeta n^{-1} + \sbeta m^{-2\alpha} \right).
		\end{align*}
	\end{remark}
	
	\begin{example}[Minimax Estimation of the Mean Function]\label{examplecycommon}
		Consider again the model given in equation \eqref{modelcy} with common sampling times, which is reproduced below for convenience.
		\begin{align*}
		y_{i}(t_{i,j}) = \beta(t_{i,j}) + \xi_{i}(t_{i,j}) + \epsilon_{i}(t_{i,j}).
		\end{align*}
		In this case, it is clear that $g(n) = 1$ satisfies assumption \ref{assumptionq=0betaint}.  Next, observe that
		\begin{align*}
		\tr\left( \left( X^\T X \right)^{-1} X^\T X \left( X^\T X \right)^{-1} \right) = n^{-1}.
		\end{align*}
		Thus, the MISE simplifies to
		\begin{align*}
		\mise = \O \left( n^{-1} + (nm)^{-1} \kbeta + m^{-2\alpha} + \kbeta^{-2\alpha} \right).
		\end{align*}
		Choosing $\kbeta \leq m - 1$ and $\kbeta \asymp m$, it follows that
		\begin{align*}
		\mise  = \O\left( n^{-1} + m^{-2\alpha} \right),
		\end{align*}
		which coincides with the minimax rate from Theorem 2.1 of \cite{cai2011}.
	\end{example}
	
	The next result is the analogue of Theorem \ref{theoremq=0hd} for the common sampling times setting.
	
	\begin{theorem}\label{theoremq=0hdcommon}
		Consider the model given in equation \eqref{modelfdcom}.  Assume \ref{assumptionp=0sgerrors}, \ref{assumptionq=0samplingtimes}, \ref{assumptionq=0designscaling}, \ref{assumptionq=0betaint}, and \ref{assumptionq=0xisup}.  For $k \leq \kbeta$ and $t > 0$, let
		\begin{align*}
		&c_{k} \defined \begin{cases}
		\sgparam_{\daleth, k}^{2} + \sqrt{2} \sum_{r = 1}^{\infty} \sgparam_{\daleth, 2rm}^{2}, & k = 1, \\
		\sgparam_{\daleth, k}^{2} + \sum_{r = 1}^{\infty} \sgparam_{\daleth, 2rm + k}^{2} + \sgparam_{\daleth, 2rm - k}^{2}, &  k = 2, 4, \dots, m - 1, \\
		\sgparam_{\daleth, k}^{2} + \sum_{r = 1}^{\infty} \sgparam_{\daleth, 2rm + k}^{2} + \sgparam_{\daleth, 2rm + 2 - k}^{2}, &  k = 3, 5, \dots, m - 1,
		\end{cases} \\
		&\lambda_{0, k} \defined \sqrt{c_{k} + m^{-1} \sgparam_{\epsilon}^{2}} \sqrt{\frac{t^2 + 2 \log(p)}{n}}.
		\end{align*}
		Suppose $\lambda_{k} \geq 2 \lambda_{0, k}$.
		\begin{enumerate}
			\item If, in addition \ref{assumptionq=0xcc} holds, then with probability at least $1 - 2 \exp(-t^{2} / 2)$,
			\begin{align*}
			\Vert X (\bethhathd_{k} - \beth_{k} - \samekh_{k}) \Vert_{2}^{2} + \lambda_{k} \Vert \bethhathd_{k} - \beth_{k} - \samekh_{k} \Vert_{1}
			\leq 4 \lambda_{k}^{2} \sbeta / \ccx^{2}.
			\end{align*}
			
			\item If, in addition \ref{assumptionq=0xarev} holds, then with probability at least $1 - 2 \exp(-t^{2} / 2)$,
			\begin{align*}
			\Vert \bethhathd_{k} - \beth_{k} - \samekh_{k} \Vert_{2}^{2} = \O(\sbeta \lambda_{k}^{2}).
			\end{align*}
			
			\item If, in addition \ref{assumptionq=0xarev} holds, $\kbeta \leq m - 1$, $\kbeta \asymp m$, and $\lambda_{k} = 2\lambda_{0, k}$, then with probability at least $1 - 2\exp(-t^{2} / 2 + \log(\kbeta))$,
			\begin{align*}
			\ise(\betahathd) = \O \left( \frac{\sbeta \log(p)}{n} + \sbeta m^{-2\alpha} \right).
			\end{align*}
		\end{enumerate}
	\end{theorem}

	\section{Two-Stage Estimation}\label{sectioncomb}
	
	So far, we have only considered the special cases where either $p=0$ or $q=0$.  In practice, it is more common to have a mixture of both time varying and time invariant covariates.  In this section, we briefly describe how to extend the theory developed in the preceding two sections to consider the general model from Section \ref{sectionintroduction}, which is reproduced below.
	\begin{align*}
	y_{i}(t_{i,j}) = x_{i}^\T \beta(t_{i,j}) + [z_{i}(t_{i,j})]^\T \gamma(t_{i,j}) + \xi_{i}(t_{i,j}) + \epsilon_{i}(t_{i,j}).
	\end{align*}
	Since $\beta(\cdot)$ is smooth and $x_{i}$ is time invariant, the product $x_{i}^\T \beta(\cdot)$ is smooth.  Thus, we may similarly consider the differencing approach from Section \ref{sectionp=0} to simultaneously remove the effect of $x_{i}^\T \beta(\cdot)$ and $\xi_{i}(\cdot)$.  Again, for simplicity, we consider the differencing given by $\Aset_{h}$.  That is, we may likewise form the differenced linear model from equation \eqref{modelplmgimelvec}
	\begin{align*}
	\Ydiff = \Psi \gimelvec + \etavec + \remgammavec,
	\end{align*}
	where
	\begin{align*}
	\remgamma_{i,j} = \phant & \sum_{k=\kgamma +1}^{\infty} [\phi_{k}(t_{i,(j+1)}) z_{i}(t_{i,(j + 1)}) - \phi_{k}(t_{i,(j)}) z_{i}(t_{i,(j)})]^\T \gimel_{k} \\
	&+ \xi_{i}(t_{i,(j+1)}) - \xi_{i}(t_{i,(j)}) + x_{i}^\T \left( \beta(t_{i,(j+1)}) - \beta(t_{i,(j)}) \right).
	\end{align*}
	Then, we may use the same estimators for $\gamma(\cdot)$ from Section \ref{sectionp=0}.  To estimate $\beta(\cdot)$, we consider the residuals after estimating $\gamma(\cdot)$.  That is, define
	\begin{align*}
	\yres_{i}(t_{i,j}) \defined y_{i}(t_{i,j}) - [z_{i}(t_{i,j})]^\T \gammahat(t_{i,j}).
	\end{align*}
	Then, we may again convert these observations to the frequency domain and use the estimators for $\beta(\cdot)$ from Section \ref{sectionq=0}.  Depending on whether the sampling times are random and independent or fixed and common, we let
	\begin{align*}
	\omega_{i,k} = x_{i}^\T \beth_{k} + \zeta_{i,k} + \rembeta_{i,k},
	\end{align*}
	or
	\begin{align*}
	\omega_{i,k} = x_{i}^\T (\beth_{k} + \samekh_{k}) + \zeta_{i,k} + \rembeta_{i,k},
	\end{align*}
	where
	\begin{align*}
	\rembeta_{i,k} \defined m_{i}^{-1} \sum_{j=1}^{m_{i}} [z_{i}(t_{i,j})]^\T \left( \gamma(t_{i,j}) - \gammahat(t_{i,j}) \right) \phi_{k}(t_{i,j}).
	\end{align*}
	The arguments are nearly identical to the preceding sections, with the only change being the change in the bias terms.  Hence, we omit the results of this section.

	\section{Confidence Bands}\label{sectioninference}
	
	In this section, we consider the problem of constructing a confidence band for a particular varying coefficient.  The problem of inference in high-dimensional varying coefficient models was first considered by \cite{chen2018}.  However, they only considered the problem of conducting inference at a prespecified time.  To the best of our knowledge, there have been no proposals for confidence bands in high-dimensions.  Since the proof technique is similar between the time invariant covariates and the time varying covariates, we only provide the details for $\beta_{1}(\cdot)$ with independent sampling times, but describe the procedure for $\gamma_{1}(\cdot)$.  For simplicity, we assume that $m_{i} = m$ for all $i = 1, \dots, n$.  Let $\sigma_{\zeta,k}^{2} \defined \var(\zeta_{i, k})$.
	
	Recall that $\beta_{1}(\cdot)$ admits a decomposition as
	\begin{align*}
	\beta_{1}(\cdot) = \betalower_{1}(\cdot) + \betaupper_{1}(\cdot),
	\end{align*}
	where $\betalower_{1}(\cdot)$ and $\betaupper_{1}(\cdot)$ are the low and high frequency components of $\beta_{1}(\cdot)$ respectively.  Under some regularity conditions, we may bound the high frequency signal by
	\begin{align*}
	\left| \sum_{k=\kbeta + 1}^{\infty} \beth_{k,1} \phi_{k}(t) \right|
	= \O\left( \kbeta^{-\alpha} \log(\kbeta) \right)
	\end{align*}
	uniformly for all $t \in (0,1)$.  Therefore, it suffices to construct a confidence band for $\beta_{1}(\cdot)$ and tune $\kbeta$ to balance with the above bias.  By the definition of $\betalower_{1}(\cdot)$, we have that
	\begin{align*}
	\betalower_{1}(\cdot) = \sum_{k=1}^{\kbeta} \beth_{k,1} \phi_{k}(\cdot).
	\end{align*}
	
	For each $k = 1, \dots, \kbeta$, a confidence interval for $\beth_{k,1}$ may be constructed using the debiased lasso.  We write $\Sigmahat \defined X^\T X / n$ and define $\Thetahat \in \R^{p \times p}$ as the relaxed inverse of $\Sigmahat$ using nodewise lasso as in \cite{vandegeer2014}.  The debiased estimator for $\beth_{k}$ is given by
	\begin{align*}
	\bethhatdb_{k} \defined \bethhathd_{k} + \Thetahat X^\T \left( \Yproj - X \bethhathd_{k} \right) / n.
	\end{align*}
	Writing $\thetahat^\T$ to denote the first row of $\Thetahat$, the debiased estimator for $\beth_{k,1}$ is
	\begin{align*}
	\bethhatdb_{k,1} \defined \bethhathd_{k,1} + \thetahat^\T X^\T \left( \Yproj - X \bethhathd_{k} \right) / n.
	\end{align*}
	By using a multiple comparisons correction procedure, we may construct simultaneous confidence intervals for $\beth_{k,1}$ for $k = 1, \dots, \kbeta$.  We may then use these simultaneous confidence intervals to extend to a confidence band for $\betalower(\cdot)$.  Let $a_{k}$ and $b_{k}$ be the lower and upper bounds for a $1 - \tau$ simultaneous confidence interval of $\beth_{k,1}$.  Then, the $1 - \tau$ lower and upper confidence bands for $\betalower_{1}(\cdot)$ will be given by
	\begin{align*}
	\lbeta(t) \defined \min_{c_{k} : \{a_{k} \leq c_{k} \leq b_{k}\}} \sum_{k=1}^{\kbeta} c_{k} \phi_{k}(t),
	&&
	\ubeta(t) \defined \max_{c_{k} : \{a_{k} \leq c_{k} \leq b_{k}\}} \sum_{k=1}^{\kbeta} c_{k} \phi_{k}(t).
	\end{align*}
	By slightly enlarging these values, we may account for the bias of $\betaupper(\cdot)$ asymptotically.  That is, for a value of $\delta > 0$, define
	\begin{align*}
	\lbeta_{\delta}(t) \defined l(t) - \delta,
	&&
	\ubeta_{\delta}(t) \defined u(t) + \delta.
	\end{align*}
	
	For $\gamma_{1}(\cdot)$, we may use a similar idea.  First, consider the simpler problem of constructing confidence intervals at an arbitrary point $\tstar \in (0, 1)$.  Let $\phivec(\tstar) = (\phi_{1}(\tstar), \dots, \phi_{\kgamma}(\tstar))^\T$ denote a loading vector, which we identify $\phivec(\tstar)$ as a vector in $\R^{\kgamma}$ as well as a vector in $\R^{q\kgamma}$.  Since, as a vector in $\R^{q\kgamma}$, $\phivec(\tstar)$ is a sparse loading vector, to construct a confidence interval interval for $\gammalower(\cdot)$, we may consider the approach of \cite{cai2017} for estimating linear functionals.  This yields a confidence interval for $\gammalower(\tstar)$.
	
	Then, consider the time grid $1 / (2\kgamma), 2 / (2\kgamma), \dots, 2\kgamma / (2\kgamma)$.  By using a multiple comparisons adjustment, we may construct simultaneous $1 - \tau$ confidence intervals at the $2\kgamma$ time points.  To extend this to a confidence band for $\gammalower(\cdot)$, consider the following lower and upper bounds:
	\begin{align*}
	\lgamma(t) \defined \min_{c_{k} : \{a_{k} \leq c_{k} \leq b_{k}\}} \sum_{k=1}^{\kgamma} c_{k} \frac{\sin(\pi(2\kgamma t - k))}{\pi(2 \kgamma t - k)},
	\\
	\ugamma(t) \defined \max_{c_{k} : \{a_{k} \leq c_{k} \leq b_{k}\}} \sum_{k=1}^{\kgamma} c_{k} \frac{\sin(\pi(2\kgamma t - k))}{\pi(2 \kgamma t - k)}.
	\end{align*}
	At first glance, this definition may seem strange.  However, our confidence band is leveraging the Nyquist-Shannon Theorem (cf. \cite{shannon1949}).  Since every low-frequency signal, with maximal frequency $\kgamma$, can be recovered by interpolation of the signal at the grid points $1 / (2\kgamma), 2 / (2\kgamma), \dots, 2\kgamma / (2\kgamma)$ using the sinc function, the above band simultaneously covers all possible interpolations that can arise given the confidence intervals for the signal value.  Then, to incorporate the bias in $\gammaupper(\cdot)$, we may again enlarge these intervals by a value of $\delta > 0$.
	
	We note that the idea of using a multiple comparisons correction to construct confidence bands is not new in the literature.  Earlier works such as \cite{knafl1985} and \cite{wu1998} use a Bonferroni adjustment at various gridpoints and interpolate over the interval by bounding the derivative.  Conversely, for $\beta(\cdot)$, we construct simultaneous confidence intervals on the Fourier coefficients, which induces simultaneous confidence intervals for all linear combinations.  Likewise, for $\gamma(\cdot)$, instead of bounding the derivatives, we exploit the Nyquist-Shannon Theorem to provide uniformity over the entire interval.

	\subsection{Assumptions}
	
	\begin{enumerate}[label=(C\arabic*)]
		\item \label{assumptionq=0inferencesparsity} The sparsity $\sbeta$ satisfies $\sbeta = o(\sqrt{n} / \log(p))$.
		
		\item \label{assumptionq=0inferencevariance} The projected error term $\zetavec$ satisfies $\sigma_{\zeta,k}^{2} \defined \var(\zeta_{i,k}) \asymp \sgparam_{\zeta,k}^{2}$.
		
		\item \label{assumptionq=0inferencevdg} The estimator $\Thetahat$ satisfies $\sqrt{n} \Vert \id_{p} - \Thetahat \Sigmahat \Vert_{\infty} = \Op( \sqrt{\log(p)} )$.
	\end{enumerate}
	
	\begin{remark}
		The first assumption \ref{assumptionq=0inferencesparsity} is the standard sparsity requirement in high-dimensional inference.  Next, \ref{assumptionq=0inferencevariance} is a technical requirement that is satisfied, for example, by the Gaussian distribution.  Sufficient conditions for \ref{assumptionq=0inferencevdg} are given in \cite{vandegeer2014}.
	\end{remark}
	
	\subsection{Main Results}
	
	\begin{theorem}\label{theoremq=0db}
		Consider the model given in equation \eqref{modelq=0}.  Assume \ref{assumptionp=0errors}, \ref{assumptionq=0samplingtimes}, \ref{assumptionq=0betaint}, \ref{assumptionq=0xisup}, and \ref{assumptionq=0inferencesparsity} -- \ref{assumptionq=0inferencevdg}.  Moreover, let the tuning parameters for lasso satisfy $\lambda_{k} \asymp \lambda_{0,k}$ and $\lambda_{k} \geq 3 \lambda_{0,k}$ for $k = 1, \dots, \kbeta$.  Similarly, letting $\nu_{j}$ for $j = 1, \dots, p$ denote the tuning parameters for nodewise lasso, assume that the parameters satisfy $\nu_{j} \asymp K_{0} \sqrt{\log(p) / n}$ and $\nu_{j} \geq 3 K_{0} \sqrt{\log(p) / n}$.
		
		For each $k = 1, \dots, \kbeta$, write
		\begin{align*}
		\sqrt{n \sigma_{\zeta,k}^{-2} } \left( \bethhatdb_{k,1} - \beth_{k,1} \right) &= W_{k} + \Delta_{k}, \\
		W_{k} & = \sigma_{\zeta, k}^{-1} \thetahat^\T X^\T \zeta_{k} / \sqrt{n}, \\
		\Delta_{k} &= \sqrt{n \sigma_{\zeta,k}^{-2}} \left( \left( I_{p} - \Thetahat \Sigmahat \right) \left( \bethhathd_{k} - \beth_{k} \right) \right)_{1} .
		\end{align*}
		Let $(V_{1}, \dots, V_{\kbeta})^\T \sim \n_{\kbeta}(\zerovec_{\kbeta}, \var((W_{1}, \dots, W_{\kbeta})^\T | X))$.
		\begin{enumerate}
			\item Letting $\mathcal{A}$ denote the class of all hyperrectangles in $\R^{\kbeta}$, then
			\begin{align*}
			\sup_{A \in \mathcal{A}} \left| \p \left( (W_{1}, \dots, W_{\kbeta})^\T \in A | X \right)  - \p \left( (V_{1}, \dots, V_{\kbeta})^\T \in A | X \right) \right| \to 0.
			\end{align*}
			
			\item Then,
			\begin{align*}
			\sup_{k=1, \dots, \kbeta} |\Delta_{k}| = \op(1).
			\end{align*}
		\end{enumerate}
	\end{theorem}
	
	\begin{remark}
		In practice, one may use the scaled lasso of \cite{sun2012} to estimate $\sigma_{\zeta,k}^{2}$, denoted by $\sigmahat_{\zeta,k}^{2}$.
	\end{remark}
	
	For concreteness, we consider simultaneous confidence intervals constructed by Bonferroni correction.  Let $z_{\tau}$ denote the $\tau$ upper quantile of a standard Gaussian distribution.  Then, the confidence intervals for $\beth_{k, 1}$ are given by $\bethhatdb_{k,1} \pm z_{\tau / (2\kbeta)} \sigmahat_{\zeta,k} / \sqrt{n}$.  Then, we have the following proposition regarding the coverage and asymptotic performance of the confidence band.
	
	\begin{proposition}\label{propositionq=0cb}
		Consider the setup of Theorem \ref{theoremq=0db} with the simultaneous confidence intervals constructed by Bonferroni correction.  Suppose that $\kbeta \asymp \min((n \log(n))^{1/(2\alpha)}$, $(nm \log(n) / g(n))^{1/(2\alpha + 2))})$ and $\delta \asymp \kbeta^{-\alpha} \log(\kbeta)$.
		
		\begin{enumerate}
			\item For $n$ sufficiently large,
			\begin{align*}
			\p\left( \forall t \in (0, 1) : \lbeta_{\delta}(t) \leq \beta_{1}(t) \leq \ubeta_{\delta}(t) \right) \geq 1 - \tau.
			\end{align*}
			
			\item For all $t \in (0, 1)$,
			\begin{align*}
			\left| \ubeta_{\delta}(t) - \lbeta_{\delta}(t) \right| \leq \max \left( (n\log(n))^{-1/2}, (nm \log(n) / g(n))^{-\alpha / (2\alpha + 2)} \log(n) \right).
			\end{align*}
		\end{enumerate}
	\end{proposition}
	
	\begin{remark}
		We note that the maximal width of the confidence band has two rates, depending on the growth rate of $m$ relative to $n$.  If $m$ is very large, then most of the error in the resultant linear model comes from the random Fourier coefficients of $\xi_{i}(\cdot)$.  Since the variance for the high-frequency components is relatively low, this encourages oversmoothing to reduce the bias.  This leads to a near parametric rate of $\sqrt{\log(n) / n}$ in the width.  Conversely, when $m$ is small, the noise is dominated by the inexact orthogonalization.  This variance accumulates as we increase the number of frequencies to be estimated, resulting in less smoothing compared to the large $m$ setting.
		
		This phenomenon is related to the two part rate for estimation, as seen in Example \ref{exampleq=0hd}.  After a certain threshold of $m$, additional observations per subject do not improve the risk in estimation since the bottleneck is in averaging the random effects.
	\end{remark}

	\section{Simulations}\label{sectionnumerics}
	
	The general model that we consider is given by
	\begin{align*}
	y_{i}(t_{i,j}) = x_{i}^\T \beta(t_{i,j}) + [z_{i}(t_{i,j})]^\T \gamma(t_{i,j}) + \xi_{i}(t_{i,j}) + \epsilon_{i}(t_{i,j}),
	\end{align*}
	Throughout, the covariates and the noise are independent and identically distributed standard Gaussian variables; that is, $x_{i} \iid \n_{p}(\zerovec_{p}, \id_{p})$, $z_{i}(t_{i,j}) \iid \n_{q}(\zerovec_{q}, \id_{q})$, and $\epsilon_{i,j}\iid \n(0, 1)$.  When generating the varying coefficients and the random effects, we use either the trigonometric basis or by the B-spline basis; the choice of the coefficients is described below.  We consider both fixed and common sampling times as well as random and independent sampling times.  When the sampling times are fixed and common (denoted ``com''), we set $t_{i,j} = j / m$ for all $i = 1, \dots, n$ and $j = 1, \dots, m$.  On the other hand, when the sampling times are random and independent (denoted ``ind''), we let $t_{i,j} \iid \U(0, 1)$.
	
	All of our tuning parameters, $\lambda$, $\kbeta$, and $\kgamma$ are chosen via five-fold cross-validation.  For each combination of parameter values that we consider, we simulate the data $200$ times and compute the high-dimensional estimator.  To evaluate the performance of our estimator, we consider three metrics:  average loss, average coverage, and average length.  Average loss is integrated squared error in estimating $\beta(\cdot)$ averaged over the trials.  Average coverage is the proportion of times the confidence band covers the true varying coefficient function $\beta_{1}(\cdot)$ with a nominal coverage of $95\%$ and average length is $\max_{t\in(0,1)}(\ubeta(t) - \lbeta(t))$ averaged over the trials.  We consider average loss for both $\gamma(\cdot)$ and $\beta(\cdot)$ while average coverage and average length are only considered for $\beta(\cdot)$.  For our simulations, we consider the two special cases of Section \ref{sectionp=0} and \ref{sectionq=0}.
	
	In the setting when $q=0$, we have
	\begin{align*}
	y_{i}(t_{i,j}) = x_{i}^\T \beta(t_{i,j}) + \xi_{i}(t_{i,j}) + \epsilon_{i}(t_{i,j}).
	\end{align*}
	When $\beta(\cdot)$ is generated from the trigonometric basis, the Fourier coefficients are given by
	\begin{align*}
	\beth_{k, l} =
	\begin{cases}
	\zeta_{k, l} ((k+1) / 2)^{-2.1} & k = 1, 3, \dots, 29 \text{ and } l = 1, \dots, \sbeta \\
	\zeta_{k, l} ((k+2) / 2)^{-2.1} & k = 2, 4, \dots, 30 \text{ and } l = 1, \dots, \sbeta \\
	0 & \text{else}
	\end{cases}
	\end{align*}
	where $\zeta_{k, l} \iid \U(-1, 1)$.  Then, we rescale $\beth_{k}$ such that $\int_{0}^{1} (x_{i}^\T \beta(t))^{2} dt = 4$ to keep the signal to noise ratio constant.  For the B-spline basis, we consider
	\begin{align*}
	\beth_{k, l} =
	\begin{cases}
	\zeta_{k, l} & k = 1, 2, 3 \text{ and } l = 1, \dots, \sbeta \\
	0 & \text{else},
	\end{cases}
	\end{align*}
	which is then rescaled to keep the signal to noise ratio constant.  Similarly, for the random effects, $\xi_{i}(\cdot)$, we rescale the coefficients so the variance is constant at one.
	
\begin{figure}[!h]
		\centering
		\includegraphics[scale=0.4]{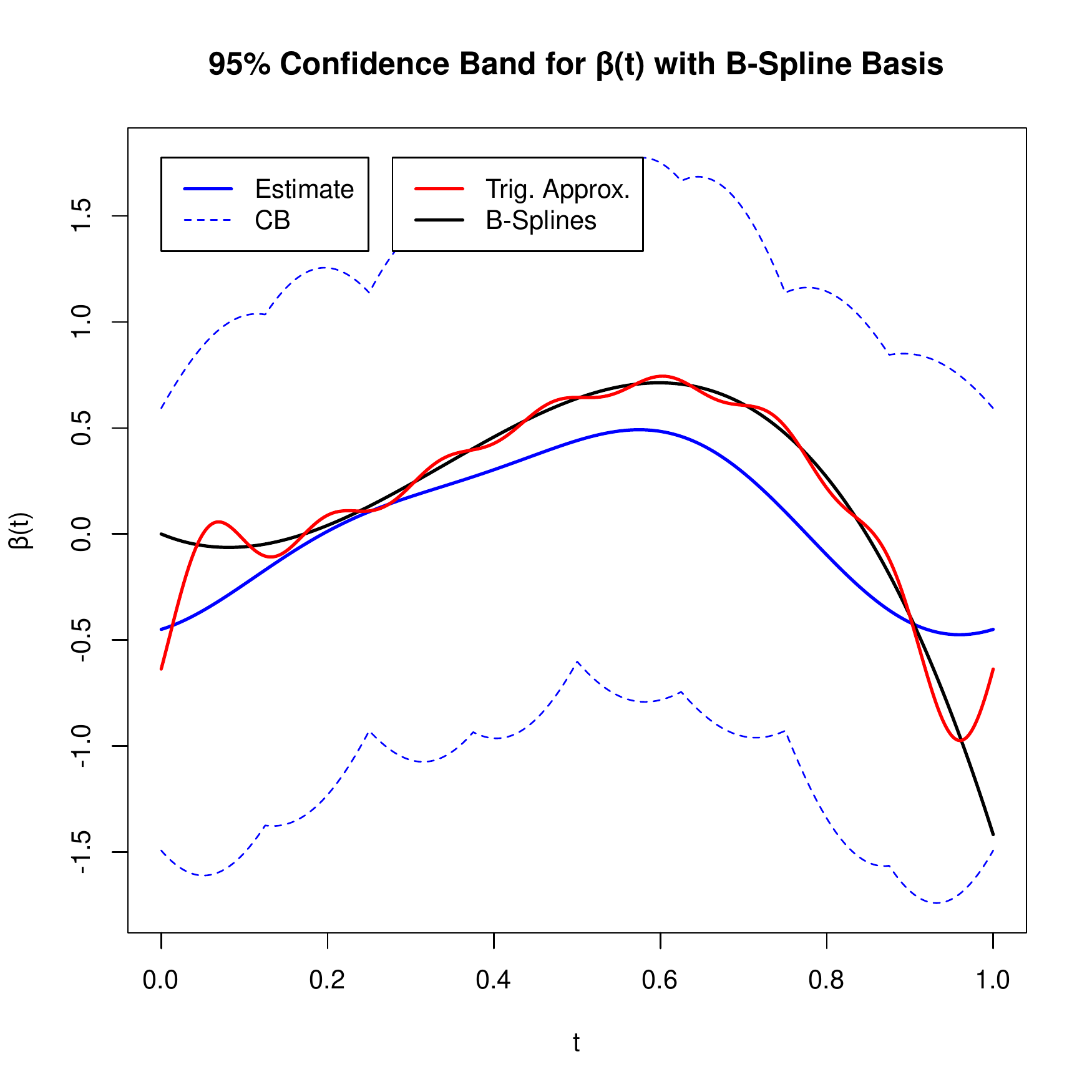}
		\includegraphics[scale=0.4]{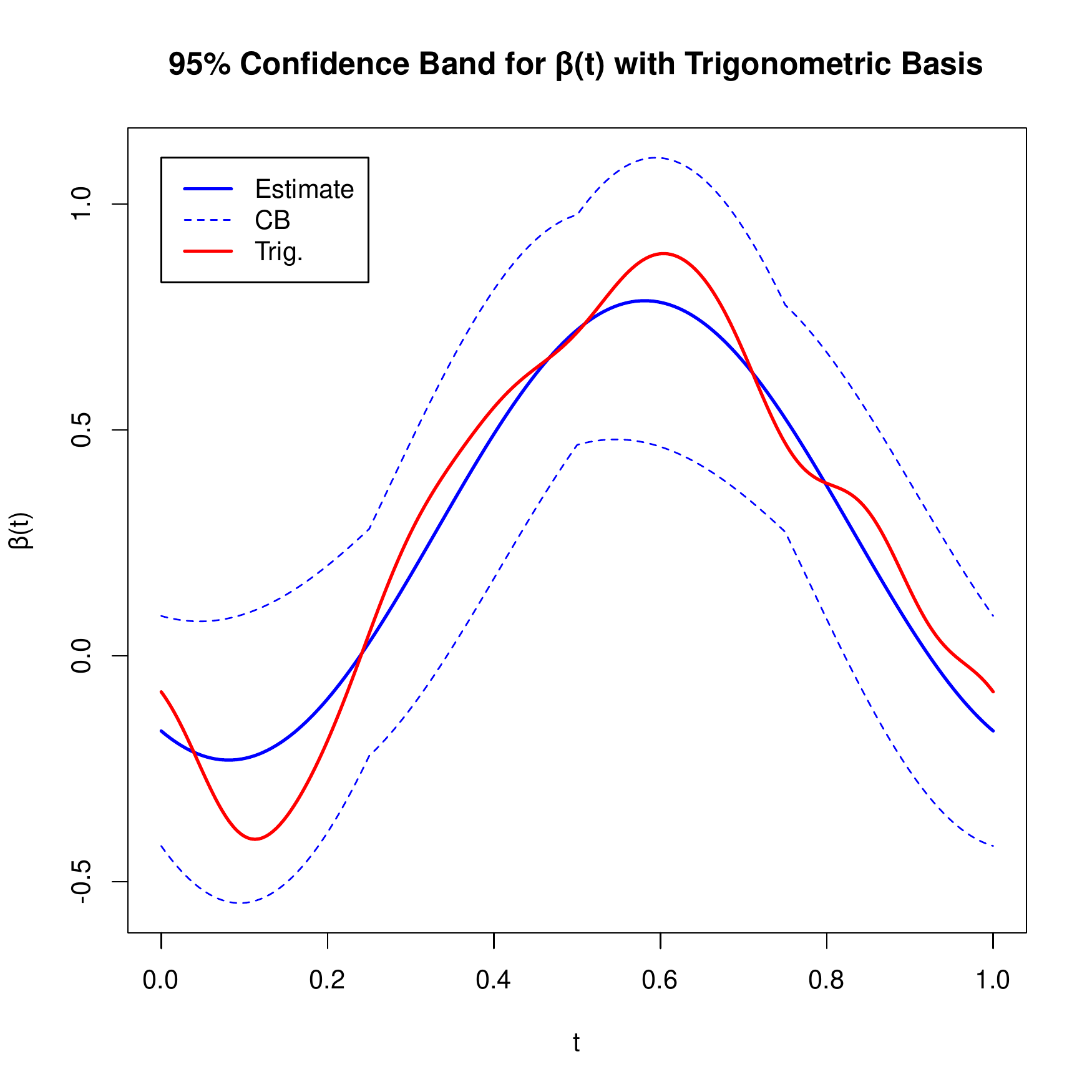}
		\caption{$95\%$ Confidence bands for $\beta_{1}(t)$ when $n = 500$, $m = 50$, $\sbeta = 15$, and $p = 1000$ with independent and random sampling times.  On the left, the data is generated using the b-spline basis.  Then, ``B-Splines'' denotes the true signal, ``Trig. Approx.'' denotes the best approximation of $\beta_{1}(t)$ using $30$ trigonometric basis function, and ``Estimate'' and ``CB'' are the estimate and confidence band from Section \ref{sectioninference} respectively.  On the right, the data is generated using the trigonometric basis.  Then, ``Trig.'' denotes the true signal and ``Estimate'' and ``CB'' are the estimate and confidence band from Section \ref{sectioninference} respectively.}
		\label{cbsplines}
	\end{figure}

	In this setting, we let $n \in \{200, 500\}$, $m \in \{25, 50, 75, 150\}$, $p \in \{500, 1000\}$, and $\sbeta \in \{15, 25\}$.  The results are presented in Tables \ref{tablebetap1000fourier} and \ref{tablebetap1000splines} in Section \ref{sectionsupplementsimulations} of the Supplement for the trigonometric and B-spline basis respectively.  For both bases, consistent with Theorems \ref{theoremq=0hd} and \ref{theoremq=0hdcommon}, as $\sbeta$ increases, the average loss increases while the loss decreases as $n$ or $m$ increase.  Surprisingly, the loss for the fixed and common sampling times seems to be better than for the random and independent sampling times for the same value of $n$ and $m$ despite the rate of convergence for the random and independent sampling times being faster.  Similarly, the confidence bands exhibit higher coverage with shorter lengths in the fixed and common sampling times as opposed to the random and independent sampling times.  We note that the confidence bands for the spline basis are significantly wider than for the trigonometric basis to account for the fact that the trigonometric basis functions periodic.  Since the spline functions are aperiodic, the confidence bands are wider to be able to cover allow coverage at both endpoints.  In Figure \ref{cbsplines}, the plot on the left shows a confidence band with the B-spline basis and the plot on the right with the trigonometric basis; note the difference in the widths of the two bands.

	In the setting when $p=0$, we have
	\begin{align*}
	y_{i}(t_{i,j}) = [z_{i}(t_{i,j})]^\T \gamma(t_{i,j}) + \xi_{i}(t_{i,j}) + \epsilon_{i}(t_{i,j}).
	\end{align*}
	Here, we set $n = 200$, $m \in \{25, 50, 75\}$, $q = 500$, and $\sgamma \in\{15, 25\}$.  The coefficients $\gimel_{k,l}$ are generated similar to $\beth_{k, l}$.  The results are presented in Table \ref{tablegammaq500} in Section \ref{sectionsupplementsimulations} of the Supplement.  We see that $\Aset_{h}$ outperforms $\Bset_{h}$ when the sampling times are independent and random whereas $\Bset_{h}$ outperforms $\Aset_{h}$ when the sampling times are common and fixed.  As $\sgamma$ increases, the estimation error increases while the estimation error decreases as $m$ increases, which is consistent with our theoretical results.  Moreover, similar to estimating $\beta(\cdot)$, our estimation error is higher for the B-spline basis as compared to the trigonometric basis.
	
	\section{Human Height Data}\label{sectionrda}
	In this section, we are interested in analyzing the average height across countries.\linebreak  The height data is freely available from the NCD Risk Factor collaboration at \linebreak https://ncdrisc.org/index.html, which includes the average height of birth cohorts aged five through nineteen over many decades for both sexes.  A detailed description of the data is provided in \cite{rodriguez2020}.  Although the data contains $95\%$ credible intervals for the average height in a country of each age group at a particular time, our response comprises solely of the point estimate.  To supplement this data, we use a variety of United Nations data on countries, including the World Health Organization (https://www.who.int/data/collections), the Human Development Reports (http://www.hdr.undp.org/en/data), and the United Nations International Children's Emergency Fund (https://data.unicef.org/).  In addition, we also use information on caloric, protein, and fat supply from Our World in Data (https://ourworldindata.org/food-supply).  For covariates that are time varying but not collected annually, we impute the values for intermediate years using an exponentially weighted moving average.   We focus on the average height of five year old boys, the youngest age in the dataset.  The model that we consider is model \eqref{modely}, which is reproduced below for convenience:
	\begin{align*}
	y_{i}(t_{i,j}) = x_{i}^\T \beta(t_{i,j}) + [z_{i}(t_{i,j})]^\T \gamma(t_{i,j}) + \xi_{i}(t_{i,j}) + \epsilon_{i}(t_{i,j}).
	\end{align*}
	
	Since some covariates are only available for a subset of countries, we only consider the $n = 98$ countries that have some measurements for all covariates.  We constrain our analysis to the years 1990 to 2019 to limit the imputation of our time-varying covariates.  Therefore, we have common sampling times, with $m = 30$.  When normalizing our sampling times to the unit interval, $t = 0$ and $t = 29 / 30$ correspond to 1990 and 2019 respectively.  Our time invariant covariates are various classifications of the countries according to the United Nations; in particular, we consider development status and geographic regions and sub-regions, with sub-regions being nested within regions.  Expanding out these groupings, we have $p = 25$ time invariant covariates.  Our time varying covariates consist of various socioeconomic factors, such as human development index, urbanization rate, gross domestic product per capita, etc.  When synthesizing the height data, \cite{finucane2014} considered the interaction between income and urbanization rate.  Thus, we include all possible two-way interactions of our time-varying covariates, leaving us with $q = 212$ covariates.  Finally, the random effects represent the unobserved heterogeneity amongst countries, which, as mentioned in the Introduction, encapsulates environmental factors.  These environmental factors are potentially correlated with our observed time-varying covariates.
	
\begin{figure}[H]
		\centering
		\includegraphics[scale=0.35]{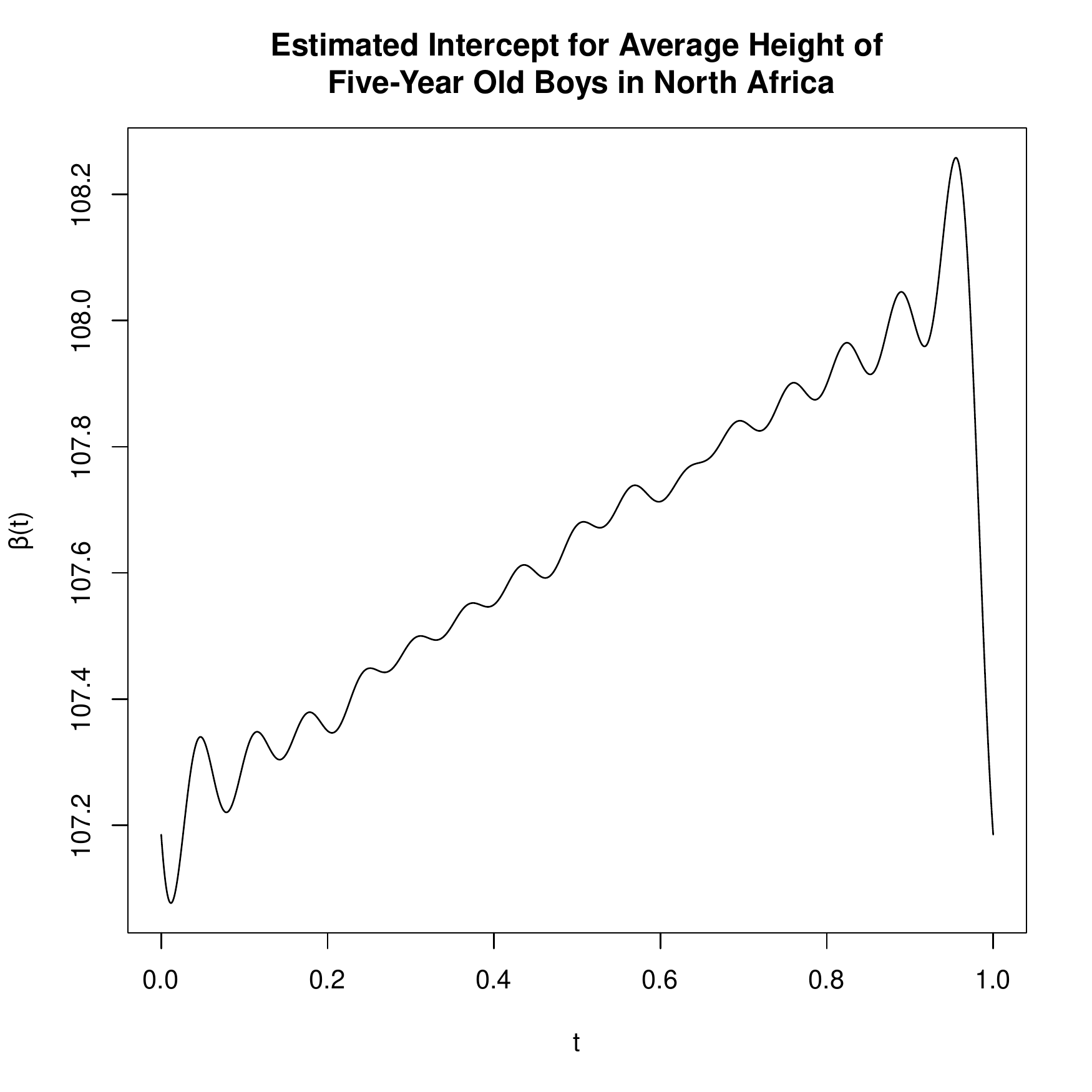}
		\caption{Estimated coefficient for the intercept, with North Africa as the reference region.}
		\label{figureheightplot}
	\end{figure}

	The goal in our data analysis is to analyze the trend in the average height, controlling for other covariates.  That is, we are interested in estimating the intercept term.  Since our time-invariant covariates consist of many geographic regions, our reference region is North Africa.  A plot of the estimate is given in Figure \ref{figureheightplot}.
	
	From the plot, we notice that the estimate of the average height in North Africa, controlling for other covariates, is almost linearly increasing, which is not specified apriori.  Note that the precipitous drop near $t = 1$ is an artifact of the estimation procedure using the trigonometric basis.  As we do not believe that the average height over the span of thirty years is periodic, our estimator $\betahat$ is estimating the best periodic approximation to the underlying non-periodic intercept function.  Thus, the interpretation of the plot is valid only away from the two endpoints.
	
	\bibliographystyle{newapa}
	\bibliography{fdaRef}

	\newpage
	

	\title{Supplement to ``High-Dimensional Varying Coefficient Models with Functional Random Effects''}
	
	\author{Michael Law \and Ya\hspace{-.1em}'\hspace{-.1em}acov Ritov}
	\date{
		University of Michigan\\
		\today
	}
	
	\maketitle
	
	\begin{abstract}
		This supplement provides additional assumptions, all of the proofs, the simulation results, and an analysis of the yeast cell cycle data.
	\end{abstract}

	\newcounter{suppsection}
	\addtocounter{suppsection}{1}
	\def\thesection{S\arabic{suppsection}}
	\section{Assumptions for Section \ref{sectionq=0}}\label{sectionsupplementassumption}
	
	In Section \ref{sectionq=0}, we use the following assumptions.
	
	\begin{enumerate}[label=(B\arabic*)]
		\item \label{assumptionq=0samplingtimes} The sampling times satisfy $t_{i,j} \iid \U(0,1)$ and are independent of the errors $\epsilon_{i}(t_{i,j})$.
		
		\item \label{assumptionq=0samplingtimesdiscrete} The number of observations per individual is the same, $m_{i} = m$ for all $i = 1, \dots, n$.  Moreover, the sampling times satisfy $t_{i,j} = j / m$ for all $i = 1, \dots, n$ and $j = 1, \dots, m$.
		
		\item \label{assumptionq=0designscaling} The columns of $X$ have squared norms that are uniformly $\O(n)$ and the entries of $X$ satisfy $\sup_{i=1, \dots, n} \Vert x_{i} \Vert_{\infty} = g(n)$ for some function $g(n)$.
		
		\item \label{assumptionq=0xcc} The design matrix $X$ satisfies the compatibility condition with compatibility constant $\ccx > 0$.
		
		\item \label{assumptionq=0xarev} The design matrix $X$ satisfies the adaptive restricted eigenvalue condition with constant $\arevx > 0$.
		
		\item \label{assumptionq=0betaint} Each coordinate of the coefficient $\beta(\cdot)$ satisfies $\beta_{k}(\cdot) \in \Wper(\alpha, R)$ for some constant $\alpha \geq 2$ and $R > 0$ with $\sbeta$ coordinates nonzero.  Moreover, the product $x_{i}^\T \beta(\cdot) \in \Wper(\alpha, \O(g(n))$ uniformly and
		\begin{align*}
		\sup_{i=1, \dots, n}\int_{0}^{1} \left( x_{i}^\T \beta(t) \right)^{2} dt = \O(g(n))
		\end{align*}
		for some function $g(n)$.
		
		\item \label{assumptionq=0xiint} The individual random effects $(\xi_{i}(\cdot))_{i=1}^{n} \subseteq \Wper(\alpha, R)$ almost surely and are independent and identically distributed.  Furthermore, for each $i = 1, \dots, n$ and $t\in (0,1)$, the function $\xi_{i}(\cdot)$ satisfies $\e \xi_{i}(t) = 0$ and
		\begin{align*}
		\e \int_{0}^{1} \xi_{i}^{2}(t) dt < \infty.
		\end{align*}
		
		\item \label{assumptionq=0xisup} The individual random effects $(\xi_{i}(\cdot))_{i=1}^{n} \subseteq \Wper(\alpha, R)$ almost surely and are independent and identically distributed.  Furthermore, for each $i = 1, \dots, n$, the Fourier coefficients $\daleth_{i,k} \sim \sg(\sgparam_{\daleth, k}^{2})$ with $\sum_{k = K}^{\infty} \sgparam_{\daleth, k}^{2} = \O(K^{-2\alpha})$.
	\end{enumerate}
	
	\addtocounter{suppsection}{1}
	\section{Some Technical Lemmata}
	
	In this section, we provide some technical lemmata regarding the trigonometric basis that are used in the proofs.  We start with a lemma regarding the fourth moments of the trigonometric basis functions.
	
	\begin{lemmaS}\label{lemmatrigbasis}
		Let $a, b, k, l \in \N$ be fixed positive integers and $(\phi_{k})_{k=1}^{\infty}$ denote the trigonometric basis as defined in Definition \ref{definitiontrigbasis}.  Then,
		\begin{align*}
		\int_{0}^{1} \phi_{a} (t) & \phi_{b}(t) \phi_{k}(t) \phi_{l}(t) dt \\
		\leq \phant & \delta_{a+b, k + l} + \delta_{a + b + k, l} + \delta_{a + b + l, k}
		+ \delta_{a + k + l, b} + \delta_{a, b + k + l} + \delta_{a + k, b + l} + \delta_{a + l, b + k} \\
		&+ \delta_{l, 1} \left( \delta_{a + b, k} + \delta_{a + k, b} + \delta_{a,  b + k} \right)
		+ \delta_{k, 1} \left( \delta_{a + b, l} + \delta_{a + l, b} + \delta_{a,  b + l} \right) \\
		&+ \delta_{b, 1} \left( \delta_{a + k, l} + \delta_{a + l, k} + \delta_{a,  k + l} \right)
		+ \delta_{a, 1} \left( \delta_{b + k, l} + \delta_{b + l, k} + \delta_{b,  k + l} \right).
		\end{align*}
	\end{lemmaS}
	\begin{proof}
		We consider a few cases:
		\begin{enumerate}
			\item $a, b, k, l > 1$ and $a, b, k, l$ are even. 	
			\begin{align*}
			\int_{0}^{1} \phi_{a} (t) & \phi_{b}(t) \phi_{k}(t) \phi_{l}(t) dt \\
			&\begin{aligned}
			= 4 \int_{0}^{1} \cos(\pi a t) \cos(\pi b t) \cos(\pi k t) \cos(\pi l t) dt
			\end{aligned} \\
			&\begin{aligned}
			= \int_{0}^{1} &(\cos(\pi (a + b) t) + \cos(\pi (a - b) t) ) \\
			&\times (\cos(\pi (k + l) t) + \cos(\pi (k - l) t)) dt
			\end{aligned} \\
			&\begin{aligned}
			=\phant & \frac{1}{2} \int_{0}^{1} \left( \cos(\pi (a + b + k + l) t)  + \cos(\pi (a + b - k - l)) \right) dt \\
			&+ \frac{1}{2} \int_{0}^{1} \left( \cos(\pi (a + b + k - l) t)  + \cos(\pi (a + b - k + l)) \right) dt \\
			&+ \frac{1}{2} \int_{0}^{1} \left( \cos(\pi (a - b + k + l) t)  + \cos(\pi (a - b - k - l)) \right) dt \\
			&+ \frac{1}{2} \int_{0}^{1} \left( \cos(\pi (a - b + k - l) t)  + \cos(\pi (a - b - k + l)) \right) dt \\
			\end{aligned} \\
			&\begin{aligned}
			= \frac{1}{2} (&\delta_{a+b, k + l} + \delta_{a + b + k, l} + \delta_{a + b + l, k} \\
			&+ \delta_{a + k + l, b} + \delta_{a, b + k + l} + \delta_{a + k, b + l} + \delta_{a + l, b + k} ).
			\end{aligned}
			\end{align*}
			
			\item $a, b, k, l > 1$, $a, b, k$ are even, and $l$ is odd. 	
			\begin{align*}
			\int_{0}^{1} \phi_{a} (t) & \phi_{b}(t) \phi_{k}(t) \phi_{l}(t) dt \\
			&\begin{aligned}
			= 4 \int_{0}^{1} \cos(\pi a t) \cos(\pi b t) \cos(\pi k t) \sin(\pi (l - 1) t) dt
			\end{aligned} \\
			&\begin{aligned}
			= \int_{0}^{1} &(\cos(\pi (a + b) t) + \cos(\pi (a - b) t) ) \\
			&\times (\sin(\pi (k + l - 1) t) - \sin(\pi (k - l + 1) t)) dt
			\end{aligned} \\
			&\begin{aligned}
			=\phant & \frac{1}{2} \int_{0}^{1} \left( \sin(\pi (a + b + k + l - 1) t)  + \sin(\pi (a + b - k - l + 1)) \right) dt \\
			&- \frac{1}{2} \int_{0}^{1} \left( \sin(\pi (a + b + k - l + 1) t)  + \sin(\pi (a + b - k + l - 1)) \right) dt \\
			&+ \frac{1}{2} \int_{0}^{1} \left( \sin(\pi (a - b + k + l - 1) t)  + \sin(\pi (a - b - k - l + 1)) \right) dt \\
			&- \frac{1}{2} \int_{0}^{1} \left( \sin(\pi (a - b + k - l + 1) t)  + \sin(\pi (a - b - k + l - 1)) \right) dt \\
			\end{aligned} \\
			&\begin{aligned}
			= 0.
			\end{aligned}
			\end{align*}
			
			\item $a, b, k, l > 1$, $a, b$ are even, and $k, l$ are odd. 	
			\begin{align*}
			\int_{0}^{1} \phi_{a} (t) & \phi_{b}(t) \phi_{k}(t) \phi_{l}(t) dt \\
			&\begin{aligned}
			= 4 \int_{0}^{1} \cos(\pi a t) \cos(\pi b t) \sin(\pi (k - 1) t) \sin(\pi (l - 1) t) dt
			\end{aligned} \\
			&\begin{aligned}
			= \int_{0}^{1} &(\cos(\pi (a + b) t) + \cos(\pi (a - b) t) ) \\
			&\times  (\cos(\pi (k - l) t) - \cos(\pi (k + l - 2) t)) dt
			\end{aligned} \\
			&\begin{aligned}
			=\phant & \frac{1}{2} \int_{0}^{1} \left( \cos(\pi (a + b + k - l) t)  + \cos(\pi (a + b - k + l)) \right) dt \\
			&- \frac{1}{2} \int_{0}^{1} \left( \cos(\pi (a + b + k + l - 2) t)  + \cos(\pi (a + b - k - l + 2)) \right) dt \\
			&+ \frac{1}{2} \int_{0}^{1} \left( \cos(\pi (a - b + k - l) t)  + \cos(\pi (a - b - k + l)) \right) dt \\
			&- \frac{1}{2} \int_{0}^{1} \left( \cos(\pi (a - b + k + l - 2) t)  + \cos(\pi (a - b - k - l + 2)) \right) dt \\
			\end{aligned} \\
			&\begin{aligned}
			\leq \frac{1}{2} (&\delta_{a + b + k, l} + \delta_{a + b + l, k} + \delta_{a + k, b + l} + \delta_{a + l, b + k}).
			\end{aligned}
			\end{align*}
			
			\item $a, b, k, l > 1$, $a$ is even, and $b, k, l$ are odd. 	
			\begin{align*}
			\int_{0}^{1} \phi_{a} (t) & \phi_{b}(t) \phi_{k}(t) \phi_{l}(t) dt \\
			&\begin{aligned}
			= 4 \int_{0}^{1} \cos(\pi a t) \sin(\pi (b - 1) t) \sin(\pi (k - 1) t) \sin(\pi (l - 1) t) dt
			\end{aligned} \\
			&\begin{aligned}
			= \int_{0}^{1} &(\sin(\pi (a + b - 1) t) - \sin(\pi (a - b + 1) t) ) \\
			&\times (\cos(\pi (k - l) t) - \cos(\pi (k + l - 2) t)) dt
			\end{aligned} \\
			&\begin{aligned}
			=\phant & \frac{1}{2} \int_{0}^{1} \left( \sin(\pi (a + b + k - l - 1) t)  + \sin(\pi (a + b - k + l - 1)) \right) dt \\
			&- \frac{1}{2} \int_{0}^{1} \left( \sin(\pi (a + b + k + l - 3) t)  + \sin(\pi (a + b - k - l + 1)) \right) dt \\
			&- \frac{1}{2} \int_{0}^{1} \left( \sin(\pi (a - b + k - l + 1) t)  + \sin(\pi (a - b - k + l + 1)) \right) dt \\
			&+ \frac{1}{2} \int_{0}^{1} \left( \sin(\pi (a - b + k + l - 1) t)  + \sin(\pi (a - b - k - l + 3)) \right) dt \\
			\end{aligned} \\
			&\begin{aligned}
			= 0.
			\end{aligned}
			\end{align*}
			
			\item $a, b, k, l > 1$ and $a, b, k, l$ are odd. 	
			\begin{align*}
			\int_{0}^{1} \phi_{a} (t) & \phi_{b}(t) \phi_{k}(t) \phi_{l}(t) dt \\
			&\begin{aligned}
			= 4 \int_{0}^{1} \sin(\pi (a - 1) t) \sin(\pi (b - 1) t) \sin(\pi (k - 1) t) \sin(\pi (l - 1) t) dt
			\end{aligned} \\
			&\begin{aligned}
			= \int_{0}^{1} &(\cos(\pi (a - b) t) - \cos(\pi (a + b - 2) t) ) \\
			&\times (\cos(\pi (k - l) t) - \cos(\pi (k + l - 2) t)) dt
			\end{aligned} \\
			&\begin{aligned}
			=\phant & \frac{1}{2} \int_{0}^{1} \left( \cos(\pi (a - b + k - l) t)  + \cos(\pi (a - b - k + l)) \right) dt \\
			&- \frac{1}{2} \int_{0}^{1} \left( \cos(\pi (a - b + k + l - 2) t)  + \cos(\pi (a - b - k - l + 2)) \right) dt \\
			&- \frac{1}{2} \int_{0}^{1} \left( \cos(\pi (a + b + k - l - 2) t)  + \cos(\pi (a + b - k + l - 2)) \right) dt \\
			&+ \frac{1}{2} \int_{0}^{1} \left( \cos(\pi (a + b + k + l - 4) t)  + \cos(\pi (a + b - k - l)) \right) dt \\
			\end{aligned} \\
			&\begin{aligned}
			\leq \frac{1}{2} (\delta_{a + b, k + l} + \delta_{a + k, b + l} + \delta_{a + l, b + k}).
			\end{aligned}
			\end{align*}
			
			\item $a, b, k > 1$, $a, b ,k$ are even, and $l = 1$.
			\begin{align*}
			\int_{0}^{1} \phi_{a} (t) & \phi_{b}(t) \phi_{k}(t) \phi_{l}(t) dt \\
			&\begin{aligned}
			= 2\sqrt{2} \int_{0}^{1} \cos(\pi a t) \cos(\pi b t) \cos(\pi k t) dt
			\end{aligned} \\
			&\begin{aligned}
			= \sqrt{2} \int_{0}^{1} &(\cos(\pi (a + b) t) + \cos(\pi (a - b) t) ) \cos(\pi k t) dt
			\end{aligned} \\
			&\begin{aligned}
			=\phant & \frac{1}{\sqrt{2}} \int_{0}^{1} \left( \cos(\pi (a + b + k) t)  + \cos(\pi (a + b - k)) \right) dt \\
			&+ \frac{1}{\sqrt{2}} \int_{0}^{1} \left( \cos(\pi (a - b + k) t)  + \cos(\pi (a - b - k)) \right) dt \\
			\end{aligned} \\
			&\begin{aligned}
			= \frac{1}{\sqrt{2}} (\delta_{a + b, k} + \delta_{a + k, b} + \delta_{a, b + k}).
			\end{aligned}
			\end{align*}
			
			\item $a, b, k > 1$, $a, b$ are even, $k$ is odd, and $l = 1$.
			
			\begin{align*}
			\int_{0}^{1} \phi_{a} (t) & \phi_{b}(t) \phi_{k}(t) \phi_{l}(t) dt \\
			&\begin{aligned}
			= 2\sqrt{2} \int_{0}^{1} \cos(\pi a t) \cos(\pi b t) \sin(\pi (k - 1) t) dt
			\end{aligned} \\
			&\begin{aligned}
			= \sqrt{2} \int_{0}^{1} &(\cos(\pi (a + b) t) + \cos(\pi (a - b) t) ) \sin(\pi (k - 1) t)  dt
			\end{aligned} \\
			&\begin{aligned}
			=\phant & \frac{1}{\sqrt{2}} \int_{0}^{1} \left( \sin(\pi (a + b + k - 1) t) - \sin(\pi (a + b - k + 1)) \right) dt \\
			&+ \frac{1}{\sqrt{2}} \int_{0}^{1} \left( \sin(\pi (a - b + k - 1) t)  - \sin(\pi (a - b - k + 1)) \right) dt \\
			\end{aligned} \\
			&\begin{aligned}
			= 0.
			\end{aligned}
			\end{align*}
			
			\item $a, b, k > 1$, $a$ is even, $b, k$ are odd, and $l = 1$.
			\begin{align*}
			\int_{0}^{1} \phi_{a} (t) & \phi_{b}(t) \phi_{k}(t) \phi_{l}(t) dt \\
			&\begin{aligned}
			= 2\sqrt{2} \int_{0}^{1} \cos(\pi a t) \sin(\pi (b - 1) t) \sin(\pi (k - 1) t) dt
			\end{aligned} \\
			&\begin{aligned}
			= \sqrt{2} \int_{0}^{1} & \cos(\pi a t) (\cos(\pi (b - k) t) - \cos(\pi (b + k - 2) t) ) dt
			\end{aligned} \\
			&\begin{aligned}
			=\phant & \frac{1}{\sqrt{2}} \int_{0}^{1} \left( \cos(\pi (a + b - k) t)  + \cos(\pi (a - b + k)) \right) dt \\
			&- \frac{1}{\sqrt{2}} \int_{0}^{1} \left( \cos(\pi (a + b + k - 2) t)  + \cos(\pi (a - b - k + 2)) \right) dt \\
			\end{aligned} \\
			&\begin{aligned}
			\leq \frac{1}{\sqrt{2}} (\delta_{a + b, k} + \delta_{a + k, b}).
			\end{aligned}
			\end{align*}
			
			\item $a, b, k > 1$, $a, b, k$ are odd, and $l = 1$.
			\begin{align*}
			\int_{0}^{1} \phi_{a} (t) & \phi_{b}(t) \phi_{k}(t) \phi_{l}(t) dt \\
			&\begin{aligned}
			= 2\sqrt{2} \int_{0}^{1} \sin(\pi (a - 1) t) \sin(\pi (b - 1) t) \sin(\pi (k - 1) t) dt
			\end{aligned} \\
			&\begin{aligned}
			= \sqrt{2} \int_{0}^{1} & \sin(\pi (a - 1) t) (\cos(\pi (b - k) t) - \cos(\pi (b + k - 2) t) ) dt
			\end{aligned} \\
			&\begin{aligned}
			=\phant & \frac{1}{\sqrt{2}} \int_{0}^{1} \left( \sin(\pi (a + b - k - 1) t)  + \sin(\pi (a - b + k - 1)) \right) dt \\
			&- \frac{1}{\sqrt{2}} \int_{0}^{1} \left( \sin(\pi (a + b + k - 3) t)  + \sin(\pi (a - b - k + 1)) \right) dt \\
			\end{aligned} \\
			&\begin{aligned}
			= 0.
			\end{aligned}
			\end{align*}
			
			\item $a, b > 1$ and $k, l = 1$. $\int_{0}^{1} \phi_{a} (t) \phi_{b}(t) \phi_{k}(t) \phi_{l}(t) dt = \delta_{a,b}$.
			
			\item $a > 1$ and $b, k, l = 1$.  $\int_{0}^{1} \phi_{a} (t) \phi_{b}(t) \phi_{k}(t) \phi_{l}(t) dt = \delta_{a, 1}$.
		\end{enumerate}
		
		Combining these cases together and considering all permutations finishes the proof.
	\end{proof}
	
	Next, we have a lemma regarding the aliasing effect in Fourier transforms on the uniform grid.
	
	\begin{lemmaS}\label{lemmaalias}
		For $k$ even and $m \in \N$,
		\begin{align*}
		m^{-1} \sum_{j=1}^{m} \exp\left( \Im \pi k j / m\right) =
		\begin{cases}
		1, & \text{if } k = cm \text{ for } c \text{ even}, \\
		0, & \text{else.}
		\end{cases}
		\end{align*}
	\end{lemmaS}
	
	\begin{proof}[Proof of Lemma \ref{lemmaalias}]
		Suppose that $k = cm$ for $c \in \mathbb{Z}$.  Then, we have that $\exp(\Im \pi k j / m) = \exp(\Im \pi c j) = (-1)^{cj}$ for all $j = 1, \dots, m$.  If $c$ is even, then
		\begin{align*}
		m^{-1} \sum_{j=1}^{m} \exp\left( \Im \pi k j / m\right) = m^{-1} \sum_{j=1}^{m} 1 = 1.
		\end{align*}
		Now, if $c$ is odd, this implies that $m$ is even since $k$ is even.  Thus,
		\begin{align*}
		m^{-1} \sum_{j=1}^{m} \exp\left( \Im \pi k j / m\right) =
		m^{-1} \sum_{j=1}^{m} (-1)^{j}
		= 0.
		\end{align*}
		Finally, suppose that $k \neq cm$ for any $c \in \mathbb{Z}$.  In this setting, we have that $\exp(\Im \pi k / m) \neq 1$ while $\exp(\Im \pi k) = 1$.  Therefore,
		\begin{align*}
		m^{-1} \sum_{j=1}^{m} \exp(\Im \pi k j / m) = m^{-1} \exp(\Im \pi k / m) \frac{1 - \exp(\Im \pi k)}{1 - \exp(\Im \pi k / m)} = 0.
		\end{align*}
		This finishes the proof.
	\end{proof}
	
	Finally, the following lemma is a refinement of Lemma 1.7 of \cite{tsybakov2008} regarding the orthogonality of the trigonometric basis on the uniform grid.
	
	\begin{lemmaS}\label{lemmatrigbasisdiscrete}
		Let $m \in \N$ and $k \leq 1, 2, \dots m - 1$.
		\begin{enumerate}
			\item If $l = 1, \dots, m - 1$, then
			\begin{align*}
			m^{-1} \sum_{j=1}^{m} \phi_{k}(j / m) \phi_{l}(j / m) = \delta_{k, l}.
			\end{align*}
			
			\item If $l = m, m+1, \dots$, then
			\begin{align*}
			m^{-1} \sum_{j=1}^{m} \phi_{k}(j / m) \phi_{l}(j / m) =
			\begin{cases}
			\sqrt{2} \indic{l / m \in 2 \mathbb{Z}}, &  k = 1, \\
			\indic{(l - k) / m \in 2\mathbb{Z}} + \indic{(l + k) / m \in 2\mathbb{Z}}, &  k = 2, 4, \dots, m - 1, \\
			\indic{(l - k) / m \in 2\mathbb{Z}} - \indic{(l + k - 2) / m \in 2\mathbb{Z}}, &  k = 3, 5, \dots, m - 1.
			\end{cases}
			\end{align*}
		\end{enumerate}
	\end{lemmaS}
	
	\begin{proof}[Proof of Lemma \ref{lemmatrigbasisdiscrete}]
		The setting where $l = 1, \dots, m - 1$ is exactly Lemma 1.7 of \cite{tsybakov2008}.  For the other setting, we consider a few separate cases.  Note that the last line in each of the following cases is a consequence of Lemma \ref{lemmaalias}.
		
		\begin{enumerate}
			\item $k = 1$ and $l$ is even.
			\begin{align*}
			m^{-1} \sum_{j=1}^{m} \phi_{k}(j / m) \phi_{l}(j / m)
			&= \sqrt{2} m^{-1} \sum_{j=1}^{m} \cos(\pi l j / m) \\
			&= \frac{1}{\sqrt{2} m} \sum_{j=1}^{m} \left( \exp \left( \Im \pi l j / m \right) + \exp \left( - \Im \pi l j / m \right) \right) \\
			&= \sqrt{2} \indic{l / m \in 2\mathbb{Z}}.
			\end{align*}
			
			\item $k = 1$ and $l$ is odd.
			\begin{align*}
			m^{-1} \sum_{j=1}^{m} \phi_{k}(j / m) \phi_{l}(j / m)
			&= \sqrt{2} m^{-1} \sum_{j=1}^{m} \sin(\pi (l - 1) j / m) \\
			&= \frac{1}{\sqrt{2} m \Im} \sum_{j=1}^{m} \left( \exp \left( \Im \pi (l - 1) j / m \right) - \exp \left( - \Im \pi (l - 1) j / m \right) \right) \\
			&= 0.
			\end{align*}
			
			\item $k, l$ are both even.
			\begin{align*}
			m^{-1} &\sum_{j=1}^{m} \phi_{k}(j / m) \phi_{l}(j / m) \\
			&\begin{aligned}
			= m^{-1} \sum_{j=1}^{m} \left( \cos(\pi (l - k) j / m) + \cos(\pi (l + k) j / m) \right)
			\end{aligned} \\
			&\begin{aligned}
			= \phant & \frac{1}{2m} \sum_{j=1}^{m} \left( \exp \left( \Im \pi (l - k) j / m \right) + \exp \left( - \Im \pi (l - k) j / m \right) \right) \\
			&+ \frac{1}{2m} \sum_{j=1}^{m} \left( \exp \left( \Im \pi (l + k) j / m \right) + \exp \left( - \Im \pi (l + k) j / m \right) \right)
			\end{aligned} \\
			&\begin{aligned}
			=  \indic{(l - k) / m \in 2\mathbb{Z}} + \indic{(l + k) / m \in 2\mathbb{Z}} .
			\end{aligned}
			\end{align*}
			
			\item $k$ is even and $l$ is odd.
			\begin{align*}
			m^{-1} &\sum_{j=1}^{m} \phi_{k}(j / m) \phi_{l}(j / m) \\
			&\begin{aligned}
			= m^{-1} \sum_{j=1}^{m} \left( \sin(\pi (l - k - 1) j / m) + \sin(\pi (l + k - 1) j / m) \right)
			\end{aligned} \\
			&\begin{aligned}
			= \phant & \frac{1}{2m \Im} \sum_{j=1}^{m} \left( \exp \left( \Im \pi (l - k - 1) j / m \right) - \exp \left( - \Im \pi (l - k - 1) j / m \right) \right) \\
			&+ \frac{1}{2m \Im} \sum_{j=1}^{m} \left( \exp \left( \Im \pi (l + k - 1) j / m \right) - \exp \left( - \Im \pi (l + k - 1) j / m \right) \right)
			\end{aligned} \\
			&\begin{aligned}
			=  0.
			\end{aligned}
			\end{align*}
			
			\item $k > 1$ is odd and $l$ is even.
			\begin{align*}
			m^{-1} &\sum_{j=1}^{m} \phi_{k}(j / m) \phi_{l}(j / m) \\
			&\begin{aligned}
			= m^{-1} \sum_{j=1}^{m} \left( \sin(\pi (l + k - 1) j / m) - \sin(\pi (l - k + 1) j / m) \right)
			\end{aligned} \\
			&\begin{aligned}
			= \phant & \frac{1}{2m \Im} \sum_{j=1}^{m} \left( \exp \left( \Im \pi (l + k - 1) j / m \right) - \exp \left( - \Im \pi (l + k - 1) j / m \right) \right) \\
			&- \frac{1}{2m \Im} \sum_{j=1}^{m} \left( \exp \left( \Im \pi (l - k + 1) j / m \right) - \exp \left( - \Im \pi (l - k + 1) j / m \right) \right)
			\end{aligned} \\
			&\begin{aligned}
			=  0.
			\end{aligned}
			\end{align*}
			
			\item $k, l > 1$ are both odd.
			\begin{align*}
			m^{-1} &\sum_{j=1}^{m} \phi_{k}(j / m) \phi_{l}(j / m) \\
			&\begin{aligned}
			= m^{-1} \sum_{j=1}^{m} \left( \cos(\pi (l - k) j / m) - \cos(\pi (l + k - 2) j / m) \right)
			\end{aligned} \\
			&\begin{aligned}
			= \phant & \frac{1}{2m} \sum_{j=1}^{m} \left( \exp \left( \Im \pi (l - k) j / m \right) + \exp \left( - \Im \pi (l - k) j / m \right) \right) \\
			&- \frac{1}{2m} \sum_{j=1}^{m} \left( \exp \left( \Im \pi (l + k - 2) j / m \right) + \exp \left( - \Im \pi (l + k - 2) j / m \right) \right)
			\end{aligned} \\
			&\begin{aligned}
			=  \indic{(l - k) / m \in 2\mathbb{Z}} - \indic{(l + k - 2) / m \in 2\mathbb{Z}} .
			\end{aligned}
			\end{align*}
		\end{enumerate}
		Combining these calculations together proves the claim.
	\end{proof}

	\addtocounter{suppsection}{1}
	\section{Proofs for Section \ref{sectionp=0}}
	
	\subsection{Proofs for Section \ref{sectionp=0samplesize}}
	
	The proof of Proposition \ref{propositionsamplesize} relies on the following two lemmata, which we state for completeness.
	
	\begin{lemmaS}[\cite{chao1972}]\label{lemmachao}
		Let $X \sim \Bin(n, p)$.  Then,
		\begin{align*}
		\e \left( \frac{1}{X + 1} \right) = \frac{1 - (1 - p)^{n+1}}{(n+1)p}.
		\end{align*}
	\end{lemmaS}
	
	\begin{lemmaS}[\cite{boland2002}]\label{lemmaboland}
		Let $Y \sim \Bin(n, p)$ and $X = \sum_{i=1}^{n} X_{i}$, where the $X_{i} \sim \Bin(1, p_{i})$ are independent.  Then, $Y$ is stochastically smaller than $X$ if and only if $p \leq (\prod_{i=1}^{n} p_{i})^{1/n}$.
	\end{lemmaS}
	
	\begin{proof}[Proof of Proposition \ref{propositionsamplesize}]
		Let $F$ denote the distribution function corresponding to $f$.  Then,
		\begin{align*}
		\t \Ntilde_{h}
		&= \t \sum_{i=1}^{n} \sum_{j=1}^{m} \indic{\left\{ \exists j' \neq j : t_{i,j'} \in (t_{i,j}, t_{i,j} + h] \right\}} \\
		&= \sum_{i=1}^{n} \sum_{j=1}^{m} \p \left( \exists j' \neq j : t_{i,j'} \in (t_{i,j}, t_{i,j} + h] \right) \\
		&= \sum_{i=1}^{n} \sum_{j=1}^{m} \int_{0}^{1} \p \left( \exists j' \neq j : t_{i,j'} \in (t, t + h] | t_{i,j} = t \right) f(t) dt \\
		&= \sum_{i=1}^{n} \sum_{j=1}^{m} \int_{0}^{1} \left( 1 - \p \left( \forall j' \neq j : t_{i,j'} \notin (t, t + h] | t_{i,j} = t \right) \right) f(t) dt \\
		&= \sum_{i=1}^{n} \sum_{j=1}^{m} \left( 1 - \int_{0}^{1} \left( 1 - F(t + h) - F(t) \right)^{m - 1} f(t) dt \right) \\
		&= \sum_{i=1}^{n} \sum_{j=1}^{m} \left( 1 - \int_{0}^{1} \left( 1 - hf(t) + o(h) \right)^{m - 1} f(t) dt \right).
		\end{align*}
		Under the setting of the first claim, note that
		\begin{align*}
		\int_{0}^{1} \left( 1 - hf(t) + o(h) \right)^{m - 1} f(t) dt = 1 - h + o(h).
		\end{align*}
		Substituting this into the previous display yields the first claim.  For the third claim note that
		For the third claim, observe that
		\begin{align*}
		\int_{0}^{1} \left( 1 - hf(t) + o(h) \right)^{m - 1} f(t) dt \asymp \exp(-mh),
		\end{align*}
		which implies that
		\begin{align*}
		\p \left( \Ntilde_{h} \neq N \right)
		\leq \sum_{i = 1}^{n} \sum_{j = 1}^{m} \int_{0}^{1} \left( 1 - hf(t) + o(h) \right)^{m - 1} f(t) dt
		\to 0.
		\end{align*}
		This proves the third claim.
		
		For the remaining case, consider a non-homogeneous Poisson process with intensity function $m f(\cdot)$ on $(0,1)$.  It is well known that the unordered arrival times have the same distribution as the $t_{i,j}$.  Let $M$ denote the corresponding point process.  It is easy to see that
		\begin{align*}
		\p \left( \exists j' \neq j : t_{i,j'} \in (t, t + h] | t_{i,j} = t \right) \asymp \p \left( M((t, t+h]) \geq 1 \right).
		\end{align*}
		But, the right hand side of the above display satisfies $\p \left( M((t, t+h]) \geq 1 \right) = m f(t) h + o(mh)$ since $M$ is a Poisson process.  Since $f$ is bounded from above and below, it follows that
		\begin{align*}
		\int_{0}^{1} \p \left( \exists j' \neq j : t_{i,j'} \in (t, t + h] | t_{i,j} = t \right) f(t) dt
		\asymp mh.
		\end{align*}
		Finally, it is left to show that $\e (N_{h} + 1)^{-1} \asymp (nm^{2}h)^{-1}$.  Without the loss of generality, assume that $m$ is even.  For $i = 1, \dots, n$ and $j = 1, \dots, m / 2$, define
		\begin{align*}
		u_{i, j} \defined t_{i,(2j+1)} - t_{i, (2j-1)}
		\end{align*}
		with the convention that $t_{i, (2m+1)} = 1$ for all $i = 1, \dots, n$.  Define the following sets of random variables:
		\begin{align*}
		&W_{i,j} \defined \indic{t_{i, (2j)} \in (t_{i, (2j-1)}, t_{i, (2j-1)} + h]}, \\
		&X_{i,j} \iid \Bin(1, \min(ch / u_{i,j}, 1)), \\
		&Y_{i} \sim \Bin(m/2, cmh),
		\end{align*}
		for $i = 1, \dots, n$ and $j = 1, \dots, m/ 2$.  Then, it follows that
		\begin{align*}
		N_{h} = \sum_{i=1}^{n} \sum_{j=1}^{m/2} W_{i,j}.
		\end{align*}
		Since $f(\cdot)$ is bounded away from zero and infinity, there exists a constant $c > 0$ such that for all $W_{i, j}$, $i = 1, \dots, n$ and $j = 1, \dots, m$, we have $\p(W_{i,j} | \{u_{\cdot, \cdot}\}) \geq \min (c h / u_{i, j}, 1)$.  Moreover, conditioned on $\{u_{\cdot, \cdot}\}$, it is easy to see that $W_{i,j}$ is independent of $W_{i', j'}$ if $(i,j) \neq (i', j')$.  Therefore, we have that $X_{i,j} \preceq W_{i,j}$, where $\preceq$ denotes stochastic ordering.  Now, by construction, $\sum_{j=1}^{m/2} u_{i,j} \leq 1$, so it follows that
		\begin{align*}
		chm \leq \left( \prod_{j=1}^{m/2} \min \left( \frac{ch}{u_{i,j}}, 1 \right) \right)^{2/m}.
		\end{align*}
		Thus, Lemma \ref{lemmaboland} implies that $Y_{i} \preceq \sum_{j=1}^{m/2} X_{i,j}$.  Combining these calculations, we see that
		\begin{align*}
		\t \left( \frac{1}{N_{h} + 1} \right)
		&= \t \left( \frac{1}{\sum_{i=1}^{n} \sum_{j=1}^{m/2} W_{i,j} + 1} \right) \\
		&\leq \t \left( \frac{1}{\sum_{i=1}^{n} \sum_{j=1}^{m/2} X_{i,j} + 1} \right) \\
		&\leq \t \left( \frac{1}{\sum_{i=1}^{n} Y_{i} + 1} \right) \\
		&= \frac{1 - (1 - chm)^{nm/2+1}}{(nm/2+1)mh/2} \\
		&\leq \frac{1}{nm^{2}h/2 + hm/2},
		\end{align*}
		where the penultimate line follows from Lemma \ref{lemmachao}.  By Jensen's inequality, we have
		\begin{align*}
		\t \left( \frac{1}{N_{h} + 1} \right)  \geq \frac{1}{\t N_{h} + 1} = \frac{1}{nm^{2}h/2 + 1}.
		\end{align*}
		Thus,
		\begin{align*}
		\frac{1}{nm^{2}h/2 + 1} \leq \t \left( \frac{1}{N_{h} + 1} \right) \leq \frac{1}{nm^{2}h/2 + mh/2}.
		\end{align*}
		Since $mh \to 0$, this finishes the proof.
	\end{proof}
	
	\subsection{Proofs for Section \ref{sectionp=0estimation}}
	
	We start by bounding the squared bias of the resultant differenced linear model from equation \eqref{modelplmgimelvec}.
	
	\begin{lemmaS}\label{lemmap=0bias}
		Consider the model given in equation \eqref{modelp=0}.  Assume \ref{assumptionp=0lipschitz} and \ref{assumptionp=0gamma}.  The bias term satisfies
		\begin{align*}
		\t \left( N^{-1} \left\Vert \remgammavec \right\Vert_{2}^{2} \right) = \O\left( \sgamma \kgamma^{-2\alpha} + \lip^{2} h^{2} \right).
		\end{align*}
	\end{lemmaS}
	
	\begin{proof}[Proof of Lemma \ref{lemmap=0bias}]
		Recall that each entry of $\remgammavec$ may be written as
		\begin{align*}
		\remgammavec_{i,j} =
		\sum_{k=\kgamma+1}^{\infty} [\phi_{k}(t_{i,(j+1)}) z_{i}(t_{i,(j + 1)}) - \phi_{k}(t_{i,(j)}) z_{i}(t_{i,(j)})]^\T \gimel_{k}
		+ \xi_{i}(t_{i,(j+1)}) - \xi_{i}(t_{i,(j)})
		\end{align*}
		for some $(i,j) \in \Aset_{h}$.  We consider the two parts separately.
		For the first term, it follows immediately from assumption \ref{assumptionp=0gamma} that
		\begin{align*}
		\t \left( \sum_{k=\kgamma+1}^{\infty} [\phi_{k}(t_{i,(j+1)}) z_{i}(t_{i,(j+1)}) - \phi_{k}(t_{i,(j)}) z_{i}(t_{i,(j)})]^\T \gimel_{k} \right)^{2} = \O(\sgamma \kgamma^{-2\alpha}).
		\end{align*}
		By assumption \ref{assumptionp=0lipschitz}, we may bound the second term by
		\begin{align*}
		\left( \xi_{i}(t_{i,(j+1)}) - \xi_{i}(t_{i,(j)}) \right)^{2} \leq \lip^{2} h^{2}.
		\end{align*}
		Combining these two bounds finishes the proof.
	\end{proof}

	Next, we prove the result for the low-dimensional setting.
	
	\begin{proof}[Proof of Proposition \ref{theoremp=0ld}]
		Indeed, note that
		\begin{align*}
		\gimelvechatld
		&= \left( \Psi^\T \Psi \right)^{-1} \Psi^\T \left( \Psi \gimelvec + \eta + \remgammavec \right) \\
		&= \gimelvec + \left( \Psi^\T \Psi \right)^{-1} \Psi^\T \eta + \left( \Psi^\T \Psi \right)^{-1} \Psi^\T \remgammavec.
		\end{align*}
		Bounding each of the two terms separately, we have for the first term that
		\begin{align*}
		\e \left\Vert \left( \Psi^\T \Psi \right)^{-1} \Psi^\T \eta \right\Vert_{2}^{2}
		= \tr \left[\left( \Psi^\T \Psi \right)^{-1} \right] \sigma_{\eta}^{2}
		= \O \left( \frac{\sgamma \kgamma}{N} \right),
		\end{align*}
		which follows from assumption \ref{assumptionp=0designld}.  For the second term, invoking Lemma \ref{lemmap=0bias} implies that
		\begin{align*}
		\t \left\Vert \left( \Psi^\T \Psi \right)^{-1} \Psi^\T \remgammavec \right\Vert_{2}^{2}
		\leq \t \left\Vert \left( \Psi^\T \Psi \right)^{-1} \Psi^\T \right\Vert_{2}^{2} \left\Vert \remgammavec \right\Vert_{2}^{2}
		= \O \left( \sgamma \kgamma^{-2\alpha} + \lip^{2} h^{2} \right),
		\end{align*}
		which finishes the proof.
	\end{proof}

	\begin{proof}[Proof of Theorem \ref{theoremp=0hd}]
		The first claim follows immediately from Theorem 6.2 of \cite{buhlmann2011}.  Then, for the second claim, applying Corollary 6.5 of \cite{buhlmann2011} yields
		\begin{align*}
		\left\Vert \gimelvechathd - \gimelvec \right\Vert_{2}^{2}
		&\leq 6 \lambda^{2} \sgamma \kgamma \left( \frac{3 \t (N^{-1} \Vert \remgammavec \Vert_{2}^{2})}{ \lambda^{2} \sgamma \kgamma} + \frac{16}{\arevpsi^{2}}\right)^{2}.
		\end{align*}
		Now, Lemma \ref{lemmap=0bias} implies that
		\begin{align*}
		\left\Vert \gimelvechathd - \gimelvec \right\Vert_{2}^{2}
		&= \O \left( \lambda^{2} \sgamma \kgamma \left( \frac{\sgamma \kgamma^{-2\alpha} + \lip^{2}h^{2}}{\lambda^{2} \sgamma \kgamma} + \arevpsi^{-2} \right)^{2} \right) ,
		\end{align*}
		which finishes the proof.
	\end{proof}

	\addtocounter{suppsection}{1}
	\section{Proofs for Section \ref{sectionq=0}}
	
	\subsection{Proofs for Section \ref{sectionq=0independent}}
	
	\begin{proof}[Proof of Proposition \ref{theoremq=0ld}]
		Throughout this proof, we consider the model given by equation \eqref{modelfdind}.  Recall from equation \eqref{equationmiseparseval} that
		\begin{align*}
		\mise (\betahatld) = \e \t \sum_{k=1}^{\kbeta} \left\Vert \bethhatld_{k} - \beth_{k} \right\Vert_{2}^{2} + \e \t \sum_{k=\kbeta+1}^{\infty} \left\Vert \beth_{k} \right\Vert_{2}^{2}.
		\end{align*}
		We consider each of the two sums separately.  Suppose temporarily that $k \leq \kbeta$.  Then, the risk in estimating $\beth_{k}$ by $\bethhatld_{k}$ is given by
		\begin{align*}
		\e \t \left\Vert \bethhatld_{k} - \beth_{k} \right\Vert_{2}^{2}
		&= \e \t \left\Vert \left( X^\T X \right)^{-1} X^\T \zetavec_{k} \right\Vert_{2}^{2} \\
		&= \tr\left( \left( X^\T X \right)^{-1} X^\T \left( \e \t \zetavec_{k} \zetavec_{k}^\T \right) X \left( X^\T X \right)^{-1} \right).
		\end{align*}
		We directly compute $\e \zetavec_{k} \zetavec_{k}^\T$.  Note that $\e \zeta_{i,k} = 0$ for all $i = 1, \dots, n$.  Then, the covariance matrix is diagonal since observations corresponding to different individuals are independent.  Therefore, it is left to compute the value of the diagonal entries.  Fix $i = 1, \dots, n$ arbitrarily.  Then,
		\begin{align*}
		\e \t \zeta_{i,k}^{2}
		\leq \phant & 3 \e \daleth_{i,k}^{2}
		+ 3 m_{i}^{-2} \e \t \left( \sum_{j=1}^{m_{i}} \epsilon_{i}(t_{i,j}) \phi_{k}(t_{i,j}) \right)^{2} \\
		&+ 3 m_{i}^{-2} \e \t \left( \sum_{j=1}^{m_{i}} \left( x_{i}^\T \beta(t_{i,j}) + \xi_{i}(t_{i,j}) \right) \phi_{k}(t_{i,j}) - x_{i}^\T \beth_{k} - \daleth_{i,k} \right)^{2}
		\end{align*}
		We bound each of the three terms separately.  By definition, we have that $\e \daleth_{i,k}^{2} = \sigma_{\daleth, k}^{2}$.  Next,
		\begin{align*}
		m_{i}^{-2} \e \t \left( \sum_{j=1}^{m_{i}} \epsilon_{i}(t_{i,j}) \phi_{k}(t_{i,j}) \right)^{2}
		= m_{i}^{-2} \sum_{j=1}^{m_{i}} \e \epsilon_{i}^{2}(t_{i,j}) \t \phi_{k}^{2}(t_{i,j})
		= m_{i}^{-1} \sigmaepsilon.
		\end{align*}
		For the last term, an expansion yields
		\begin{align*}
		\e \t &\left( \sum_{j=1}^{m_{i}} \left( x_{i}^\T \beta(t_{i,j}) + \xi_{i}(t_{i,j}) \right) \phi_{k}(t_{i,j}) - x_{i}^\T \beth_{k} - \daleth_{i,k}\right)^{2} \\
		&\begin{aligned}
		= \e \t \left( \sum_{a=1}^{\infty} \left(x_{i}^\T \beth_{a} + \daleth_{i,a} \right) \left( \sum_{j=1}^{m_{i}} \phi_{a}(t_{i,j}) \phi_{k}(t_{i,j}) - \delta_{a,k} \right) \right)^{2}
		\end{aligned} \\
		&\begin{aligned}
		= \e \sum_{a=1}^{\infty} \sum_{b=1}^{\infty}
		&\left(x_{i}^\T \beth_{a} + \daleth_{i,a} \right)
		\left(x_{i}^\T \beth_{b} + \daleth_{i,b} \right) \\
		&\times \t \left( \sum_{j=1}^{m_{i}} \phi_{a}(t_{i,j}) \phi_{k}(t_{i,j}) - \delta_{a,k} \right)
		\left( \sum_{j=1}^{m_{i}} \phi_{b}(t_{i,j}) \phi_{k}(t_{i,j}) - \delta_{b,k} \right)
		\end{aligned}
		\end{align*}
		Applying Lemma \ref{lemmatrigbasis} shows that
		\begin{align*}
		\t &\left( \sum_{j=1}^{m_{i}} \phi_{a}(t_{i,j}) \phi_{k}(t_{i,j}) - \delta_{a,k} \right)
		\left( \sum_{j=1}^{m_{i}} \phi_{b}(t_{i,j}) \phi_{k}(t_{i,j}) - \delta_{b,k} \right) \\
		&\leq m_{i} ( \delta_{a+b, 2k} + \delta_{a + 2k, b} + \delta_{a, b + 2k} + 2\delta_{a, b} + 2 \delta_{a + 1, b} + 2 \delta_{a,  b + 1} + \delta_{a, 2k} \delta_{b,1} + \delta_{b,  2k} \delta_{a, 1}).
		\end{align*}
		By substitution, we have the following bound
		\begin{align*}
		\e \t &\left( \sum_{j=1}^{m_{i}} \left( x_{i}^\T \beta(t_{i,j}) + \xi_{i}(t_{i,j}) \right) \phi_{k}(t_{i,j}) - x_{i}^\T \beth_{k} - \daleth_{i,k}\right)^{2} \\
		&\begin{aligned}
		\leq \phant & m_{i} \e \sum_{a=1}^{2k - 1}
		\left|x_{i}^\T \beth_{a} + \daleth_{i,a} \right|
		\left|x_{i}^\T \beth_{2k - a} + \daleth_{i,2k - a} \right| \\
		&\begin{aligned}
		&+ 2 m_{i} \e \sum_{a=1}^{\infty}
		\left|x_{i}^\T \beth_{a} + \daleth_{i,a} \right|
		\left|x_{i}^\T \beth_{2k + a} + \daleth_{i,2k + a} \right| \\
		&+ 2 m_{i} \e \sum_{a=1}^{\infty}
		\left( x_{i}^\T \beth_{a} + \daleth_{i,a} \right)^{2} \\
		&+ 4 m_{i} \e \sum_{a=1}^{\infty}
		\left|x_{i}^\T \beth_{a} + \daleth_{i,a} \right|
		\left|x_{i}^\T \beth_{a + 1} + \daleth_{i,a + 1} \right| \\
		&+ 2 m_{i} \e \left|x_{i}^\T \beth_{1} + \daleth_{i,1} \right|
		\left|x_{i}^\T \beth_{2k} + \daleth_{i,2k} \right|.
		\end{aligned}
		\end{aligned}
		\end{align*}
		From assumptions \ref{assumptionq=0betaint} and \ref{assumptionq=0xiint} in conjunction with Parseval's Theorem, it follows that
		\begin{align*}
		2 m_{i} \e \sum_{a=1}^{\infty}
		\left( x_{i}^\T \beth_{a} + \daleth_{i,a} \right)^{2} = \O(m_{i}g(n)).
		\end{align*}
		Now, using the inequality $2uv \leq u^{2} + v^{2}$ and the above, we have that
		\begin{align*}
		2 m_{i} &\e \sum_{a=1}^{\infty} \left|x_{i}^\T \beth_{a} + \daleth_{i,a} \right| \left|x_{i}^\T \beth_{2k + a} + \daleth_{i,2k + a} \right| \\
		&\leq m_{i} \e \sum_{a=1}^{\infty} \left( x_{i}^\T \beth_{a} + \daleth_{i,a} \right)^{2} + m_{i} \e \sum_{a=1}^{\infty} \left( x_{i}^\T \beth_{2k + a} + \daleth_{i,2k + a} \right)^{2} \\
		&= \O(m_{i} g(n)).
		\end{align*}
		Similarly,
		\begin{align*}
		&m_{i} \e \sum_{a=1}^{2k - 1} \left|x_{i}^\T \beth_{a} + \daleth_{i,a} \right| \left|x_{i}^\T \beth_{2k - a} + \daleth_{i,2k - a} \right|
		= \O(m_{i} g(n)), \\
		&4 m_{i} \e \sum_{a=1}^{\infty} \left|x_{i}^\T \beth_{a} + \daleth_{i,a} \right| \left|x_{i}^\T \beth_{a + 1} + \daleth_{i,a + 1} \right|
		= \O(m_{i} g(n)).
		\end{align*}
		Thus, combining all the results yields
		\begin{align*}
		\e \t \left\Vert \bethhatld_{k} - \beth_{k} \right\Vert_{2}^{2}
		= \tr\left( \left( X^\T X \right)^{-1} X^\T \left( \sigma_{\daleth, k}^{2} I_{n} + \O(g(n)) M^{-1} \right) X \left( X^\T X \right)^{-1} \right),
		\end{align*}
		where $M \in \R^{n \times n}$ is a diagonal matrix whose $i$'th entry is $m_{i}$.  Therefore,
		\begin{align*}
		\sum_{k=1}^{\kbeta} &\left\Vert \bethhatld_{k} - \beth_{k} \right\Vert_{2}^{2} = \tr\left( \left( X^\T X \right)^{-1} X^\T \left(
		\O\left(1\right) I_{n}
		+ \O(g(n) \kbeta)M^{-1} \right)
		X \left( X^\T X \right)^{-1} \right),
		\end{align*}
		since $\sum_{k=1}^{\kbeta} \sigma_{\daleth,k}^{2} = \O(1)$.
		Recalling that
		\begin{align*}
		\sum_{k=\kbeta+1}^{\infty} \left\Vert \beth_{k} \right\Vert_{2}^{2} = \O(\sbeta \kbeta^{-2\alpha})
		\end{align*}
		finishes the proof.
	\end{proof}

	\begin{lemmaS}\label{propositionsg}
		Consider the model from equation \eqref{modelfdind}.  Assume \ref{assumptionp=0sgerrors}, \ref{assumptionq=0betaint}, and \ref{assumptionq=0xisup}.  Then, $\zeta_{i,k} \sim \sg(\sgparam_{\zeta,i,k}^{2})$ with respect to the joint probability measure on $\epsilon_{i}(t_{i,j})$, $\xi_{i}(\cdot)$, and $t_{i,j}$, where $\sgparam_{\zeta,i,k}^{2} = \O(\sgparam_{\daleth, k}^{2} + g(n) m_{i}^{-1})$.
		
	\end{lemmaS}
	\begin{proof}[Proof of Lemma \ref{propositionsg}]
		Consider first the model from equation \eqref{modelfdind}.  We partition $\zeta_{i,k}$ into four terms and show each term is sub-Gaussian.
		\begin{align*}
		\zeta_{i,k}
		= \phant & \underbrace{\daleth_{i,k}}_{(I)}
		+ \underbrace{m_{i}^{-1} \sum_{j=1}^{m_{i}} \epsilon_{i}(t_{i,j}) \phi_{k}(t_{i,j})}_{(II)}
		+ \underbrace{m_{i}^{-1} \sum_{j=1}^{m_{i}} x_{i}^\T \left( \beta(t_{i,j}) \phi_{k}(t_{i,j}) - \beth_{k} \right)}_{(III)} \\
		&+ \underbrace{m_{i}^{-1} \sum_{j=1}^{m_{i}} \left( \xi_{i}(t_{i,j}) \phi_{k}(t_{i,j}) - \daleth_{i,k} \right)}_{(IV)}.
		\end{align*}
		For the first term, by assumption \ref{assumptionq=0xisup}, $\daleth_{i,k} \sim \sg(\sgparam_{\daleth,k}^{2})$.  Next, we have, for any fixed $\lambda > 0$,
		\begin{align*}
		\e \t \exp \left( \lambda m_{i}^{-1} \sum_{j=1}^{m_{i}} \epsilon_{i,j} \phi_{k}(t_{i,j}) \right)
		&= \prod_{j=1}^{m_{i}} \e \t \exp\left( \lambda m_{i}^{-1} \epsilon_{i,j} \phi_{k}(t_{i,j}) \right) \\
		&= \prod_{j=1}^{m_{i}} \t \exp\left( \frac{\sgparam_{\epsilon}^{2} \lambda^{2} \phi_{k}^{2}(t_{i,j})}{2 m_{i}^{2}} \right) \\
		&\leq \exp \left( \frac{2\sgparam_{\epsilon}^{2} m_{i}^{-1} \lambda^{2}}{2} \right).
		\end{align*}
		In the second equality, we have used assumption \ref{assumptionp=0sgerrors}.  Thus, $(II) \sim \sg(2 \sgparam_{\epsilon}^{2}m_{i}^{-1})$.  Then, for the third term, we see that for $\lambda > 0$,
		\begin{align*}
		\e \t & \exp \left( \lambda m_{i}^{-1} \sum_{j=1}^{m_{i}} x_{i}^\T \left( \beta(t_{i,j}) \phi_{k}(t_{i,j}) - \beth_{k} \right) \right) \\
		&\leq \e \prod_{j=1}^{m_{i}} \t \exp \left( \lambda m_{i}^{-1} x_{i}^\T \left( \beta(t_{i,j}) \phi_{k}(t_{i,j}) - \beth_{k} \right) \right) \\
		&\leq \e \prod_{j=1}^{m_{i}} \exp \left( \frac{\O(g(n)) m_{i}^{-2} \lambda^{2}}{2} \right) \\
		&\leq \exp \left( \frac{\O(g(n) m_{i}^{-1}) \lambda_{2}^{2}}{2} \right).
		\end{align*}
		Hence, $m_{i}^{-1} \sum_{j=1}^{m_{i}} x_{i}^\T \left( \beta(t_{i,j}) \phi_{k}(t_{i,j}) - \beth_{k} \right) \sim \sg(\O(g(n) m_{i}^{-1})$.  Finally, from Lemma 1.8 of \cite{tsybakov2008}, assumption \ref{assumptionq=0xisup} implies that $\left(\xi_{i}(\cdot)\right)_{i=1}^{n}$ are uniformly bounded by a constant, which we temporarily denote by $c > 0$.  Then, by an analogous argument as above, it follows that
		\begin{align*}
		m_{i}^{-1} \sum_{j=1}^{m_{i}} \left( \xi_{i}(t_{i,j}) \phi_{k}(t_{i,j}) - \daleth_{i,k} \right) \sim \sg(c^{2} m_{i}^{-1}).
		\end{align*}
		Combining these results finishes the proof.
	\end{proof}

	\begin{proof}[Proof of Theorem \ref{theoremq=0hd}]
		We proceed by modifying the standard lasso arguments to account for the different sub-Gaussian parameters of the noise term.  From the Basic Inequality (Lemma 6.1 of \cite{buhlmann2011}), it follows that
		\begin{align*}
		n^{-1} \left\Vert X\left( \bethhat_{k} - \beth_{k} \right) \right\Vert_{2}^{2} + \lambda_{k} \left\Vert \bethhat_{k} \right\Vert_{1}
		\leq 2n^{-1} \zetavec_{k}^\T X \left( \bethhat_{k} - \beth_{k} \right) + \lambda_{k} \left\Vert \beth_{k} \right\Vert_{1}.
		\end{align*}
		To bound the first term on the right hand side, we similarly apply an $\ell_{1}-\ell_{\infty}$ bound to obtain
		\begin{align*}
		2n^{-1} \left| \zetavec_{k}^\T X \left( \bethhat_{k} - \beth_{k} \right) \right|
		\leq 2n^{-1} \max_{j=1,\dots, p} \left| \zetavec_{k}^\T X_{j} \right| \left\Vert \bethhat_{k} - \beth_{k} \right\Vert_{1}.
		\end{align*}
		Define the set $\Tset_{k}$ as
		\begin{align*}
		\Tset_{k} \defined \left\{ 2n^{-1} \max_{j=1,\dots, p} \left| \zetavec_{k}^\T X_{j} \right| \leq \lambda_{0,k} \right\}.
		\end{align*}
		Then, for any value of $r > 0$,
		\begin{align*}
		\p \left( \Tset_{k}^\C \right)
		&\leq 2 \sum_{j=1}^{p} \p \left( 2n^{-1} \zetavec^\T_{k} X_{j} > \lambda_{0,k} \right) \\
		&\leq 2 \sum_{j=1}^{p} \p \left( \exp \left( r \zetavec_{k}^\T X_{j} \right) > \exp \left( \frac{rn\lambda_{0,k}}{2} \right) \right) \\
		&\leq 2 \sum_{j=1}^{p} \exp \left( - \frac{rn \lambda_{0,k}}{2} \right) \e \exp \left( r \zetavec_{k}^\T X_{j} \right) \\
		&\leq 2 \sum_{j=1}^{p} \exp \left( - \frac{rn \lambda_{0,k}}{2} \right) \exp \left( \frac{r^{2}}{2} \sum_{i=1}^{n} \sgparam_{\zeta,i,k}^{2} x_{i,j}^{2} \right)
		\end{align*}
		Since the value of $r > 0$ was arbitrary, setting
		\begin{align*}
		r = \left( \sum_{i=1}^{n} \sgparam_{\eta, i, k}^{2} x_{i,j}^{2} \right)^{-1} \frac{n \lambda_{0,k}}{2}
		\end{align*}
		yields the bound
		\begin{align*}
		\p \left( \Tset^\C \right)
		\leq 2p \exp \left( -\left( t^{2}/2 + \log(p) \right) \right) \leq 2 \exp\left( -t^{2}/2 \right).
		\end{align*}
		Thus, we may restrict our attention to the event $\Tset_{k}$.  The desired results follow by applying Theorem 6.1 and Corollary 6.5 of \cite{buhlmann2011} respectively.
	\end{proof}

	\subsection{Proofs for Section \ref{sectionq=0common}}
	
	\begin{lemmaS}\label{lemmaq=0biascommon}
		Assume \ref{assumptionq=0samplingtimesdiscrete} and \ref{assumptionq=0betaint}.  Then,
		\begin{align*}
		\sum_{k = 1}^{m - 1} \Vert \samekh_{k} \Vert_{2}^{2} = \O(\sbeta m^{-2\alpha}).
		\end{align*}
	\end{lemmaS}
	
	\begin{proof}[Proof of Lemma \ref{lemmaq=0biascommon}]
		Indeed, by the second half of Lemma \ref{lemmatrigbasisdiscrete}, it follows that
		\begin{align*}
		\samekh_{k} =
		\begin{cases}
		\sqrt{2} \sum_{r=1}^{\infty} \beth_{2rm}, & \text{if } k = 1, \\
		\sum_{r=1}^{\infty} \beth_{2rm + k} + \beth_{2rm - k}, & \text{if } k = 2, 4, \dots, m - 1, \\
		\sum_{r=1}^{\infty} \beth_{2rm + k} - \beth_{2rm + 2 - k}, & \text{if } k = 3, 5, \dots, m - 1.
		\end{cases}
		\end{align*}
		Define the following sequence of constants $(a_{k})_{k=1}^{\infty}$ from \cite{tsybakov2008}
		\begin{align*}
		a_{k} =
		\begin{cases}
		k^{\alpha}, & \text{for even } k, \\
		(k - 1)^{\alpha}, & \text{for odd } k.
		\end{cases}
		\end{align*}
		Since $\alpha > 1/2$, let $c > 0$ be a constant such that $\sum_{r=1}^{\infty} r^{-2\alpha} \leq c$.  Now, for $k = 1$, it follows that
		\begin{align*}
		\Vert \samekh_{k} \Vert_{2}^{2}
		\leq 2\left( \sum_{r=1}^{\infty} a_{2rm}^{2} \Vert \beth_{2rm} \Vert_{2}^{2} \right) \left( \sum_{r=1}^{\infty} a_{2rm}^{-2} \right)
		\leq 2 cm^{-2\alpha} \left( \sum_{r=1}^{\infty} a_{2rm}^{2} \Vert \beth_{2rm} \Vert_{2}^{2} \right).
		\end{align*}
		Similarly, for $k = 2, \dots, m - 1$ and $k = 3, \dots, m -1$, we have
		\begin{align*}
		\Vert \samekh_{k} \Vert_{2}^{2}
		\leq 2cm^{-2\alpha} \sum_{r = 1}^{\infty} \left( a_{2rm + k}^{2} \Vert \beth_{2rm + k} \Vert_{2}^{2} + a_{2rm - k}^{2} \Vert \beth_{2rm - k} \Vert_{2}^{2} \right)
		\end{align*}
		and
		\begin{align*}
		\Vert \samekh_{k} \Vert_{2}^{2}
		\leq 2cm^{-2\alpha} \sum_{r = 1}^{\infty} \left( a_{2rm + k}^{2} \Vert \beth_{2rm + k} \Vert_{2}^{2} + a_{2rm + 2 - k}^{2} \Vert \beth_{2rm + 2 - k} \Vert_{2}^{2} \right)
		\end{align*}
		respectively.  Thus, combining the above calculations yields
		\begin{align*}
		\sum_{k=1}^{m - 1} \Vert \samekh_{k} \Vert_{2}^{2}
		\leq 2 cm^{-2\alpha} \sum_{r=m}^{\infty} a_{r}^{2} \Vert \beth_{r} \Vert_{2}^{2}
		= \O(\sbeta m^{-2\alpha}),
		\end{align*}
		which finishes the proof.
	\end{proof}
	
	\begin{proof}[Proof of Proposition \ref{theoremq=0ldcommon}]
		Note that MISE is given by
		\begin{align*}
		\mise(\betahatld)
		&= \sum_{k=1}^{\kbeta} \e \Vert \bethvechatld_{k} - \bethvec_{k} \Vert_{2}^{2} + \sum_{k = \kbeta + 1}^{\infty} \Vert \bethvec_{k} \Vert_{2}^{2} \\
		&\leq 2\sum_{k=1}^{\kbeta} \left( \Vert \samekh_{k} \Vert_{2}^{2} + \e \Vert (X^\T X)^{-1} X^\T \zetavec_{k} \Vert_{2}^{2} \right) + \sum_{k = \kbeta + 1}^{\infty} \Vert \bethvec_{k} \Vert_{2}^{2}.
		\end{align*}
		For the variance term, we have that
		\begin{align*}
		\e \Vert (X^\T X)^{-1} X^\T \zetavec_{k} \Vert_{2}^{2}
		= \tr \left( (X^\T X)^{-1} X^\T \e (\zetavec_{k} \zetavec_{k}^\T) X (X^\T X)^{-1} \right).
		\end{align*}
		By independence, for $1 \leq i < j \leq n$, it follows that
		\begin{align*}
		\e (\etavec_{k} \etavec_{k}^\T)_{i,j} = 0.
		\end{align*}
		Thus, $\e (\etavec_{k} \etavec_{k}^\T)$ is a diagonal matrix.  For $i = 1, \dots, n$ and $k = 1$,
		\begin{align*}
		\e (\etavec_{k} \etavec_{k}^\T)_{i, i}
		&= \e \left( \daleth_{i, k} + m^{-1} \sum_{j=1}^{m} \phi_{k}(t_{i,j}) \epsilon_{i}(t_{i,j}) + \sqrt{2} \sum_{r=1}^{\infty} \daleth_{i, 2rm} \right)^{2} \\
		&= \e \left( \daleth_{i, k} + \sqrt{2} \sum_{r = 1}^{\infty} \daleth_{i, 2rm} \right)^{2} + \sigmaepsilon m^{-1} \\
		&\leq 2c \left( a_{k}^{2} \e \daleth_{i, k}^{2} + \sum_{r = 1}^{\infty} a_{2rm}^{2} \e \daleth_{i, 2rm}^{2} \right) + \sigmaepsilon m^{-1}.
		\end{align*}	
		By performing similar calculations when $k = 2, \dots, m - 1$ and $k = 3, \dots, m - 1$, we have that
		\begin{align*}
		\e (\etavec_{k} \etavec_{k}^\T)_{i, i}\leq \begin{cases}
		2c \left( a_{k}^{2} \e \daleth_{i, k}^{2} + \sum_{r = 1}^{\infty} a_{2rm}^{2} \e \daleth_{i, 2rm}^{2} \right) + \sigmaepsilon m^{-1}, & k = 1, \\
		2c \sum_{r = 1}^{\infty}\left(  a_{2rm + k}^{2} \e \daleth_{i, 2rm + k}^{2} +  a_{2rm - k}^{2} \e \daleth_{i, 2rm}^{2} \right) + \sigmaepsilon m^{-1}, & k \text{ even}\\
		2c  \sum_{r = 1}^{\infty}\left(  a_{2rm + k}^{2} \e \daleth_{i, 2rm + k}^{2} +  a_{2rm + 2 - k}^{2} \e \daleth_{i, 2rm + 2 - k}^{2} \right) + \sigmaepsilon m^{-1}, & \text{otherwise.}
		\end{cases}
		\end{align*}	
		Since $\xi_{i}(\cdot) \in \Wper(\alpha, R)$ almost surely by assumption \ref{assumptionq=0xiint}, it follows from Proposition 1.14 of \cite{tsybakov2008} that
		\begin{align*}
		\sum_{r = 1}^{\infty} a_{k}^{2} \e \daleth_{i, r}^{2} = \O(1).
		\end{align*}
		Thus,
		\begin{align*}
		\sum_{k = 1}^{\kbeta} \e (\etavec_{k} \etavec_{k}^\T)_{i, i}
		= \O(1 + \kbeta m^{-1}).
		\end{align*}
		Hence,
		\begin{align*}
		\sum_{k = 1}^{\kbeta} \e \Vert (X^\T X)^{-1} X^\T \etavec_{k} \Vert_{2}^{2} = \O(1 + \kbeta m^{-1}) \tr [(X^\T X)^{-1}].
		\end{align*}
		Invoking Lemma \ref{lemmaq=0biascommon} and recalling that
		\begin{align*}
		\sum_{k = \kbeta + 1}^{\infty} \Vert \bethvec_{k} \Vert_{2}^{2}
		= \O(\sbeta \kbeta^{-2\alpha})
		\end{align*}
		finishes the proof.
	\end{proof}
	
	\begin{proof}[Proof of Theorem \ref{theoremq=0hdcommon}]
		We start by showing that $\zeta_{i,k}$ is sub-Gaussian with parameter $c_{k} + m^{-1} \sgparam_{\epsilon}^{2}$.  By the second half of Lemma \ref{lemmatrigbasisdiscrete}, we have that
		\begin{align*}
		\zeta_{i,k} =
		\begin{cases}
		\daleth_{i, k} + m^{-1} \sum_{j=1}^{m} \phi_{k}(t_{i,j}) \epsilon_{i}(t_{i,j}) + \sqrt{2} \sum_{r=1}^{\infty} \daleth_{i, 2rm}, & k = 1, \\
		\daleth_{i, k} + m^{-1} \sum_{j=1}^{m} \phi_{k}(t_{i,j}) \epsilon_{i}(t_{i,j}) + \sum_{r=1}^{\infty} \daleth_{i, 2rm + k} + \daleth_{i, 2rm - k}, & k = 2, 4, \dots, m - 1, \\
		\daleth_{i, k} + m^{-1} \sum_{j=1}^{m} \phi_{k}(t_{i,j}) \epsilon_{i}(t_{i,j}) + \sum_{r=1}^{\infty} \daleth_{i, 2rm + k} - \daleth_{i, 2rm + 2 - k}, &  k = 3, 5, \dots, m - 1.
		\end{cases}
		\end{align*}
		By the first half of Lemma \ref{lemmatrigbasisdiscrete}, it is easy to see that
		\begin{align*}
		m^{-1} \sum_{j=1}^{m} \phi_{k}(t_{i,j}) \epsilon_{i}(t_{i,j}) \sim \sg(m^{-1} \sgparam_{\epsilon}^{2}).
		\end{align*}
		Then, by the triangle inequality and assumption \ref{assumptionq=0xisup}, it follows that
		\begin{align*}
		\zeta_{i,k} \sim \sg (c_{k} + m^{-1} \sgparam_{\epsilon}^{2}).
		\end{align*}
		Next, the ISE can be bounded by
		\begin{align*}
		\ise(\betahathd)
		&= \sum_{k = 1}^{\kbeta} \Vert \bethhathd_{k} - \beth_{k} \Vert_{2}^{2} + \sum_{k = \kbeta + 1}^{\infty} \Vert \beth_{k} \Vert_{2}^{2} \\
		&\leq 2 \sum_{k = 1}^{\kbeta} \Vert \bethhathd_{k} - \beth_{k} - \samekh_{k} \Vert_{2}^{2} + 2 \sum_{k = 1}^{\kbeta} \Vert \samekh_{k} \Vert_{2}^{2} + \sum_{k = \kbeta + 1}^{\infty} \Vert \beth_{k} \Vert_{2}^{2} \\
		&\leq 2 \sum_{k = 1}^{\kbeta} \Vert \bethhathd_{k} - \beth_{k} - \samekh_{k} \Vert_{2}^{2} + \O(\sbeta m^{-2\alpha} + \sbeta \kbeta^{-2\alpha}),
		\end{align*}
		where we have used Lemma \ref{lemmaq=0biascommon} in the last line.  Now, by Corollary 6.5 of \cite{buhlmann2011}, it follows that, with probability at least $1 - 2\exp(-t^{2} / 2)$,
		\begin{align*}
		\Vert \bethhathd_{k} - \beth_{k} - \samekh_{k} \Vert_{2}^{2}
		= \O(\sbeta \lambda_{k}^{2}).
		\end{align*}
		By assumption \ref{assumptionq=0xisup},
		\begin{align*}
		\sum_{k = 1}^{\kbeta} c_{k} = \O(1).
		\end{align*}
		Assuming $\lambda_{k} = 2 \lambda_{0, k}$, we have that
		\begin{align*}
		\ise(\betahathd) = \O \left( \frac{\sbeta \log(p)}{n} + \frac{\sbeta \kbeta \log(p)}{mn} + \sbeta m^{-2\alpha} + \sbeta \kbeta^{-2\alpha} \right).
		\end{align*}
		Choosing $\kbeta \asymp (mn / \log(p))^{1 / (2\alpha + 1)}$,
		\begin{align*}
		\ise(\betahathd) = \O \left( \frac{\sbeta \log(p)}{n} + \sbeta \left( \frac{\log(p)}{mn} \right)^{2\alpha / (2\alpha + 1)} + \sbeta m^{-2\alpha} \right).
		\end{align*}
	\end{proof}

	\addtocounter{suppsection}{1}
	\section{Proofs for Section \ref{sectioninference}}
	
	\begin{proof}[Proof of Theorem \ref{theoremq=0db}]
		Indeed, for each $k = 1, \dots, \kbeta$, we may rewrite $\bethhatdb_{k}$ as
		\begin{align*}
		\sqrt{n \sigma_{\zeta,k}^{-2} } \left( \bethhatdb_{k} - \beth_{k} \right)
		= \phant & \sqrt{n \sigma_{\zeta,k}^{-2} } \left( \bethhat_{k} - \beth_{k} + \Thetahat X^\T \left( X \beth_{k} - X\bethhathd_{k} + \zetavec_{k} \right) / n \right) \\
		= \phant &
		\sqrt{n \sigma_{\zeta,k}^{-2} } \left( \id_{p} - \Thetahat \Sigmahat \right) \left( \bethhathd_{k} - \beth_{k} \right)
		\\ &+
		\sqrt{n \sigma_{\zeta,k}^{-2} } \Thetahat X^\T \zetavec_{k} / n.
		\end{align*}
		Now, since $\bethhatdb_{k,1}$ is the first entry of $\bethhatdb_{k}$, we set $W_{k}$ to be the first entry of $\sqrt{n \sigma_{\zeta,k}^{-2} } \Thetahat X^\T \zetavec_{k} / n$ and $\Delta_{k}$ to be the first entry of $\sqrt{n \sigma_{\zeta,k}^{-2} } \left( \id_{p} - \Thetahat \Sigmahat \right) \left( \bethhathd_{k} - \beth_{k} \right)$.  The first claim follows by Proposition 2.1 of \cite{chernozhukov2017}.  Then, an $\ell_{1}-\ell_{\infty}$ bound yields
		\begin{align*}
		\sup_{k=1, \dots, \kbeta} |\Delta_{k}| \leq \sqrt{n} \left\Vert \id_{p} - \Thetahat \Sigmahat \right\Vert_{\infty} \sup_{k=1, \dots, \kbeta} \left\Vert \bethhathd_{k} - \beth_{k} \right\Vert_{1}.
		\end{align*}
		For the other term, Theorem \ref{theoremq=0hd} implies, with probability at least $1 - 2\exp(-\log^{2}(p) / 2 + \log(\kbeta)) \to 1$, that
		\begin{align*}
		\sup_{k=1, \dots, \kbeta} \left\Vert \bethhathd_{k} - \beth_{k} \right\Vert_{1} \leq 4 \sbeta / \ccx^{2} \sup_{k=1,\dots, \kbeta} \lambda_{k} = \O\left( \sbeta \sqrt{\log(p) / n} \right).
		\end{align*}
		Combining these bounds, we see that
		\begin{align*}
		\sup_{k=1, \dots, \kbeta} |\Delta_{k}| = \Op(\sbeta \log(p) / \sqrt{n}) = \op(1),
		\end{align*}
		which finishes the proof.
	\end{proof}
	
	\begin{proof}[Proof of Proposition \ref{propositionq=0cb}]
		Recall the decomposition for $\beta_{1}(\cdot)$ as
		\begin{align*}
		\beta_{1}(\cdot) = \betalower_{1}(\cdot) + \betaupper_{1}(\cdot).
		\end{align*}
		Note that the event
		\begin{align*}
		\{\forall t \in (0, 1) &: l(t) \leq \betalower_{1}(t) \leq u(t) \}
		\cap
		\left\{ \forall t \in (0, 1) : |\betaupper_{1}(t)| \leq \delta \right\} \\
		&\subseteq \left\{ \forall t \in (0, 1) : l_{\delta}(t) \leq \betalower_{1}(t) \leq u_{\delta}(t) \right\}.
		\end{align*}
		
		For the high-frequency signal, we have under assumption \ref{assumptionq=0betaint} that
		\begin{align*}
		\left| \betaupper(t) \right|
		= \O\left( \kbeta^{-\alpha} \log(\kbeta) \right)
		\end{align*}
		from Chapter 1.21 of \cite{jackson1941} and Section 87 of \cite{achieser1992}.  Since $\delta \asymp \kbeta^{-\alpha} \log(\kbeta)$, for $n$ sufficiently large, the event $\left\{ \forall t \in (0, 1) : |\betaupper_{1}(t)| \leq \delta \right\}$ occurs with probability one.  Moreover, for $n$ sufficiently large, Theorem \ref{theoremq=0db} implies that
		\begin{align*}
		\p\left( \forall t \in (0, 1) : l(t) \leq \betalower_{1}(t) \leq u(t) \right)
		\geq 1 - \tau.
		\end{align*}
		This proves the first claim.  For the second claim, note that
		\begin{align*}
		\sup_{t\in (0,1)} \left| u_{\delta}(t) - l_{\delta}(t) \right|
		&= \sup_{t \in (0,1)} \sum_{k=1}^{\kbeta} \left( b_{k} - a_{k} + 2\delta \right) |\phi_{k}(t)| \\
		&\leq \sqrt{2} \left( \sum_{k=1}^{\kbeta} \left( b_{k} - a_{k} \right) + 2 \kbeta \delta \right).
		\end{align*}
		Now, since $z_{\tau / \kbeta} = \O(\sqrt{\log(\kbeta)})$, assumption \ref{assumptionq=0inferencevariance} implies that
		\begin{align*}
		\sum_{k=1}^{\kbeta} \left( b_{k} - a_{k} \right)
		= \Op \left( \sqrt{\log(\kbeta) / n} + \kbeta \sqrt{g(n) \log(\kbeta) / (nm)} \right).
		\end{align*}
		Moreover, note that $\log(\kbeta) \leq \log(n)$.  Thus, we have that
		\begin{align*}
		\sup_{t\in (0,1)} &\left| u_{\delta}(t) - l_{\delta}(t) \right| = \Op\left( \sqrt{\log(n) / n} + \kbeta \sqrt{g(n) \log(n) / (nm)} + \kbeta^{-\alpha} \log(n) \right).
		\end{align*}
		Substituting the choice of $\kbeta$ finishes the proof.
	\end{proof}

	\addtocounter{suppsection}{1}
	\section{Yeast Cell Cycle Data}
	In this section, we apply our methodology to analyze transcription factors affecting the cell cycle of yeast.  Versions of this data was previously analyzed in the high-dimensional varying coefficients framework by \cite{wei2011} and \cite{bai2019}, to which we refer the interested reader for a more detailed description of the data.  For our analysis, we use the data from \cite{bai2019}, which consists of $n = 47$ genes and $p = 96$ transcription factors.  The response $y$ is the mRNA level, measured seven minutes apart for 119 minutes, yielding $m = 18$ time points for each gene.  Since the time points are evenly spaced, we set $t_{i,j} = j / m$.  There are no time varying covariates in this data.  The model that we consider is
	\begin{align*}
	y_{i}(t_{i,j}) = x_{i}^\T \beta(t_{i,j}) + \xi_{i}(t_{i,j}) + \epsilon_{i}(t_{i,j}).
	\end{align*}
	This is analogous to the model fit by \cite{bai2019}, who instead combine the noise $\xi_{i}(t_{i,j}) + \epsilon_{i}(t_{i,j})$ and assume an AR(1) covariance structure.  In Figure \ref{figureyeastconfidencebands}, we provide marginal confidence bands for two selected transcription factors: ABF1 and MAC1.
	
	\addtocounter{figure}{2}
	\begin{figure}[H]
		\centering
		\includegraphics[scale=0.40]{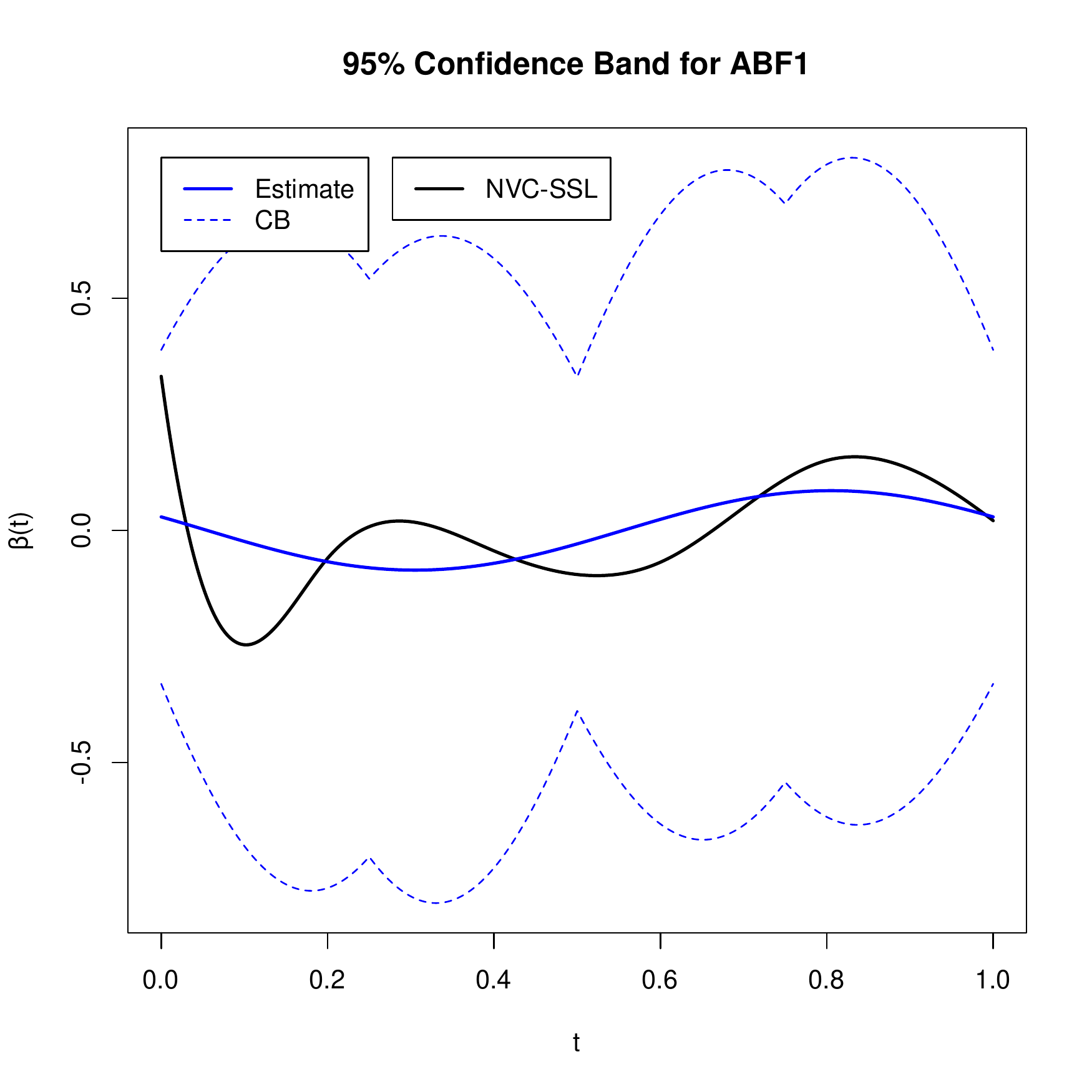}
		\includegraphics[scale=0.40]{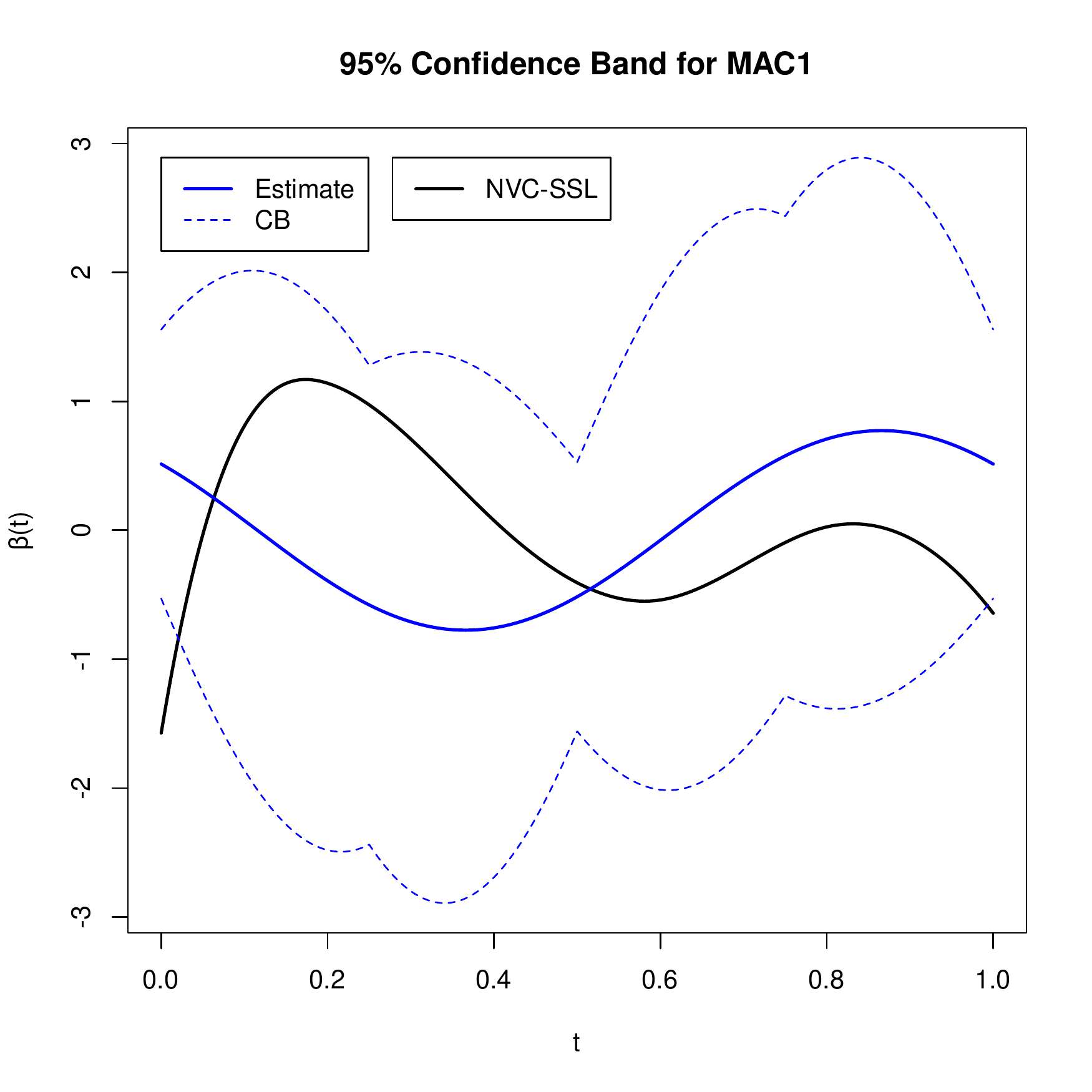}
		\caption{Marginal confidence bands for ABF1 and MAC1.  ``Estimate'' and ``CB'' are the estimate and confidence band from Section \ref{sectioninference} respectively while ``NVC-SSL'' is the methodology from \protect\cite{bai2019}.}
		\label{figureyeastconfidencebands}
	\end{figure}
	
	From the plot for ABF1, the estimated function seems to follow the general shape of the estimate from \cite{bai2019}, with both estimates entirely contained within the confidence bands.  For MAC1, we observe that the estimated curves differ for the two estimation procedures, but the confidence band contains the majority of the estimated curve by NVC-SSL, with the discrepancies at the two endpoints.  We remark that the performance of our confidence bands may be anomalous at the boundary points since the theory in Section \ref{sectionq=0} assumes that the varying coefficients are periodic Sobolev functions.  In the setting where coefficients are not periodic, convergence still holds on the interior of the interval.
	
	\addtocounter{suppsection}{1}
	\section{Additional Simulation Results}\label{sectionsupplementsimulations}
	
	In this section, we provide the results of the simulations from Section \ref{sectionnumerics}.
	
\begin{table}
\caption{\label{tablebetap1000fourier} Simulations for $\beta(\cdot)$ with Trigonometric Basis}
\centering
\fbox{
\begin{tabular}{llrrrrrrrr}
 & $\sbeta$ & 15 & 25 & 15 & 25 & 15 & 25 & 15 & 25 \\ 
   & $n$ & 200 & 200 & 500 & 500 & 200 & 200 & 500 & 500 \\ 
   & $t_{i,j}$ & ind & ind & ind & ind & com & com & com & com \\ 
   \hline
 & $m = 25$ & 0.395 & 0.582 & 0.195 & 0.246 & 0.150 & 0.252 & 0.065 & 0.086 \\ 
  Average & $m = 50$ & 0.242 & 0.368 & 0.118 & 0.149 & 0.091 & 0.153 & 0.038 & 0.051 \\ 
  Loss & $m = 75$ & 0.182 & 0.276 & 0.089 & 0.114 & 0.068 & 0.114 & 0.028 & 0.038 \\ 
   & $m = 150$ & 0.116 & 0.180 & 0.054 & 0.068 & 0.042 & 0.071 & 0.017 & 0.023 \\ 
 & $m = 25$ & 0.685 & 0.625 & 0.810 & 0.815 & 0.965 & 0.945 & 0.970 & 0.975 \\ 
  Average & $m = 50$ & 0.830 & 0.780 & 0.940 & 0.940 & 0.975 & 0.985 & 0.990 & 0.985 \\ 
  Coverage & $m = 75$ & 0.885 & 0.850 & 0.950 & 0.945 & 0.970 & 0.960 & 0.970 & 0.985 \\ 
   & $m = 150$ & 0.925 & 0.920 & 0.945 & 0.965 & 0.990 & 0.995 & 0.985 & 0.980 \\ 
 & $m = 25$ & 1.368 & 1.333 & 1.012 & 0.933 & 1.098 & 1.024 & 0.928 & 0.846 \\ 
  Average & $m = 50$ & 1.266 & 1.124 & 0.963 & 0.844 & 1.042 & 0.969 & 0.831 & 0.795 \\ 
  Length & $m = 75$ & 1.159 & 1.001 & 0.968 & 0.873 & 0.986 & 0.887 & 0.770 & 0.694 \\ 
   & $m = 150$ & 1.019 & 0.869 & 0.841 & 0.734 & 0.868 & 0.789 & 0.681 & 0.644 \\ 
\end{tabular}}
\end{table}

\begin{table}

\caption{\label{tablebetap1000splines}Simulations for $\beta(\cdot)$ with B-Spline Basis}
\centering
\fbox{\begin{tabular}{llrrrrrrrr}
 & $\sbeta$ & 15 & 25 & 15 & 25 & 15 & 25 & 15 & 25 \\ 
   & $n$ & 200 & 200 & 500 & 500 & 200 & 200 & 500 & 500 \\ 
   & $t_{i,j}$ & ind & ind & ind & ind & com & com & com & com \\ 
   \hline
 & $m = 25$ & 1.492 & 2.140 & 0.738 & 0.940 & 1.031 & 1.645 & 0.443 & 0.570 \\ 
  Average & $m = 50$ & 1.257 & 1.883 & 0.584 & 0.764 & 0.970 & 1.551 & 0.383 & 0.506 \\ 
  Loss & $m = 75$ & 1.158 & 1.755 & 0.522 & 0.682 & 0.930 & 1.490 & 0.372 & 0.489 \\ 
   & $m = 150$ & 1.031 & 1.625 & 0.432 & 0.580 & 0.889 & 1.451 & 0.337 & 0.452 \\ 
 & $m = 25$ & 0.870 & 0.755 & 0.945 & 0.945 & 1.000 & 1.000 & 1.000 & 1.000 \\ 
  Average & $m = 50$ & 0.930 & 0.910 & 0.980 & 0.985 & 1.000 & 1.000 & 1.000 & 1.000 \\ 
  Coverage & $m = 75$ & 0.990 & 0.970 & 0.990 & 1.000 & 1.000 & 1.000 & 1.000 & 1.000 \\ 
   & $m = 150$ & 1.000 & 1.000 & 1.000 & 1.000 & 1.000 & 1.000 & 1.000 & 1.000 \\ 
 & $m = 25$ & 3.755 & 3.296 & 3.449 & 3.050 & 3.751 & 3.809 & 2.870 & 2.868 \\ 
  Average & $m = 50$ & 3.923 & 3.645 & 3.709 & 3.406 & 4.013 & 4.089 & 2.915 & 2.904 \\ 
  Length & $m = 75$ & 4.587 & 4.208 & 3.583 & 3.447 & 4.407 & 4.554 & 3.130 & 3.038 \\ 
   & $m = 150$ & 5.155 & 4.863 & 4.096 & 4.015 & 4.709 & 4.871 & 3.309 & 3.296 \\ 
\end{tabular}}
\end{table}

\begin{table}

\caption{\label{tablegammaq500}Simulations for $\gamma(\cdot)$} 
\centering
\fbox{\begin{tabular}{llrrrrrrrr}
 & $\sgamma$ & 15 & 25 & 15 & 25 & 15 & 25 & 15 & 25 \\ 
   & $t_{i,j}$ & ind & ind & com & com & ind & ind & com & com \\ 
   & diff type & A & A & A & A & B & B & B & B \\ 
   \hline
Trigonometric & $m = 25$ & 1.4526 & 2.1085 & 0.3633 & 0.4806 & 0.7348 & 0.9096 & 1.3686 & 1.5289 \\ 
  Splines & $m = 25$ & 3.0488 & 4.0465 & 0.9057 & 1.1462 & 1.2152 & 1.4874 & 1.5288 & 1.7016 \\ 
  Trigonometric & $m = 50$ & 0.4231 & 0.5647 & 0.2007 & 0.2598 & 0.3992 & 0.4844 & 1.9385 & 1.9891 \\ 
  Splines & $m = 50$ & 1.0675 & 1.3570 & 0.5743 & 0.7048 & 0.6556 & 0.7819 & 2.0637 & 2.1226 \\ 
  Trigonometric & $m = 75$ & 0.2453 & 0.3167 & 0.1493 & 0.1899 & 0.3692 & 0.4238 & 2.0780 & 2.1017 \\ 
  Splines & $m = 75$ & 0.6807 & 0.8494 & 0.4437 & 0.5311 & 0.5615 & 0.6430 & 2.1772 & 2.2064 \\ 
\end{tabular}}
\end{table}

\end{document}